\documentclass[11pt,a4paper]{article}
\usepackage{amsmath, amsfonts, amssymb,}
\usepackage{a4wide}
\usepackage{parskip}
\usepackage{enumitem}
\usepackage[pagewise]{lineno}\nolinenumbers
\def\proof{\noindent{\textbf{Proof. }}}
\def\QED{\hfill {$\square$}\goodbreak \medskip}

\newtheorem{Theorem}{Theorem}[section]
\newtheorem{Lemma}[Theorem]{Lemma}
\newtheorem{Proposition}[Theorem]{Proposition}

\numberwithin{equation}{section}

\begin{document}
	\date{}
	{\vspace{0.01in}
	\title{Multiple positive solutions for degenerate Kirchhoff equations with singular and Choquard nonlinearity}
	\author{{\bf S. Rawat\footnote{email: {\tt sushmita.rawat1994@gmail.com}}  \;and  \bf K. Sreenadh\footnote{	e-mail: {\tt sreenadh@maths.iitd.ac.in}}} \\ Department of Mathematics,\\ Indian Institute of Technology Delhi,\\Hauz Khaz, New Delhi-110016, India. }
	\maketitle
	\begin{abstract}
	In this paper we study the existence, multiplicity and regularity of positive weak solutions for the following Kirchhoff-Choquard problem:
	\begin{equation*}
	 \begin{array}{cc}
			\displaystyle  M\left( \iint\limits_{\mathbb{R}^{2N}} \frac{|u(x)-u(y)|^2}{|x-y|^{N+2s}}\,dxdy\right) (-\Delta)^s u = \frac{\lambda}{u^\gamma} + \left( \int\limits_{\Omega}  \frac{|u(y)|^{2^{*}_{\mu ,s}}}{|x-y|^ \mu}\, dy\right)  |u|^{2^{*}_{\mu ,s}-2}u \;\text{in} \; \Omega,\\
			\quad \quad u = 0\quad \text{in} \; \mathbb{R}^{N}\backslash\Omega,
		\end{array}
	\end{equation*}
where $\Omega$ is open bounded domain of $\mathbb{R}^{N}$ with $C^2$ boundary, $N > 2s$ and $s \in (0,1)$. $M$ models Kirchhoff-type coefficient in particular, the degenerate case where Kirchhoff coefficient M is zero at zero. $(-\Delta)^s$ is fractional Laplace operator, $\lambda > 0$ is a real parameter, $\gamma \in (0,1)$ and $2^{*}_{\mu ,s} = \frac{2N-\mu}{N-2s}$ is the critical exponent in the sense of Hardy-Littlewood-Sobolev inequality.
We prove that each positive weak solution is bounded and satisfy H\"older regularity of order $s$. Furthermore, using the variational methods and truncation arguments we prove the existence of two positive solutions. 
\medskip

\noindent \textbf{Key words:} fractional Laplacian, Hardy-Littlewood-Sobolev critical exponent, Kirchhoff equation.

\medskip

\noindent \textit{2010 Mathematics Subject Classification: 35A15, 35J60, 35J20.}

\end{abstract}
	
	\section{Introduction}
The purpose of this article is to study a class of Kirchhoff-type equations for fractional Laplacian with singular and Choquard term. We consider the problem
\begin{equation*}    (P_\lambda)\; \left\{  \begin{array}{cc}
		\displaystyle \left( \iint\limits_{\mathbb{R}^{2N}} \frac{|u(x)-u(y)|^2}{|x-y|^{N+2s}}\,dxdy \right) ^{\theta - 1}(-\Delta)^s u = \frac{\lambda}{u^\gamma} + \left( \int\limits_{\Omega} \frac{|u(y)|^{2^{*}_{\mu ,s}}}{|x-y|^ \mu}\, dy\right)  |u|^{2^{*}_{\mu ,s}-2}u, \text{in} \; \Omega,\\
		\;u>0\; \text{in}\; \Omega,\;  u = 0\quad \text{in} \; \mathbb{R}^{N}\backslash\Omega,
	\end{array} \right.
\end{equation*}
where $\Omega$ is open bounded domain of $\mathbb{R}^{N}$ having $C^2$ boundary, $N > 2s$ with $s, \gamma \in (0,1)$, $\lambda > 0$ is a real parameter and $\theta \in [1, 2^{*}_{\mu ,s})$,  $2^{*}_{\mu ,s} = \frac{2N-\mu}{N-2s}$. Here, $(-\Delta)^s$ is fractional Laplace operator defined as,
\begin{equation*}
	(-\Delta)^su(x) = \int\limits_{\mathbb{R}^N}\dfrac{2u(x)-u(x+y)-u(x-y)}{|y|^{N+2s}}\,dy, \quad x\in \mathbb{R}^N,
\end{equation*}
A lot of research has been done on non local operators because of its many applications including it being the infinitesimal generators of the $L\acute{e}vy$ stable diffusion process \cite{Levy} and the interested readers can refer \cite{Hitchhiker's, fractionalS}  and the references therein.

On the other hand results by Crandall, Rabinowitz and Tartar \cite{crandall} were the pioneer of the work on semi-linear elliptic with singular nonlinearities. Later, various researchers including Haitao and Hirano \cite{Haitao, Hirano} studied the singular problems. In \cite{Hirano} authors   further improved the result by removing the restriction on the regularity of boundary and in \cite{Hirano Brezis} authors obtained Brezis Nirenberg type results for singular elliptic problem, with approach based on non smooth analysis. In \cite{tuhinaFS} authors used variational methods to show the existence and multiplicity of positive solutions of the following problem with respect to the parameter $\lambda$.
\begin{equation*}
	 (-\Delta)^s u = \lambda a(x)u^{-q} + u^{2^{*}_{s}-1},\quad u > 0 \;\text{in} \; \Omega,\quad u = 0\; \text{in} \; \mathbb{R}^{N} \backslash \Omega.
\end{equation*}
where $\Omega$ is a bounded domain in $\mathbb{R}^n$ with smooth boundary $\partial\Omega$, $n > 2s$, $s \in (0, 1)$, $\lambda >
0$, $0 < q \leq 1$, $\theta \leq a(x) \in L^{\infty}(\Omega)$, for some $\theta > 0$ and $2_s^{*} = \frac{2n}{n-2s}$.

Due to vast application of Choquard type non-linearity in physical models extensive study has been done in the existence and uniqueness results. One of the first applications was given by Pekar in the framework of quantum theory \cite{pekar} and Lieb \cite{choqlieb} used it in approximation of Hartree-Fock theory. For detailed understanding one can refer \cite{Moroz4, Moroz3, Moroz2, Moroz1} and references therein. In the local case $(s=1)$, Gao and Yang \cite{M.yang} analyzed Brezis-Nirenberg type existence results for Choquard type critical nonlinearity . In \cite{tuhinaFC}, authors studied these results in the following  non-local problem:
\begin{equation}\label{ts}
	(-\Delta)^s u = \lambda u^q + \left( \int\limits_{\Omega}  \frac{|u(y)|^{2^{*}_{\mu, s}}}{|x-y|^ \mu}\, dy\right)  |u|^{2^{*}_{\mu,s}-2}u \quad\text{in} \; \Omega,\; u = 0 \; \text{on} \;\mathbb{R}^{N} \backslash \Omega, 
\end{equation}
where $q=1$, $\Omega$ is a bounded domain in $\mathbb{R}^n$ with Lipschitz boundary, $\lambda$ is a real parameter, ${s\in(0, 1)}$, $n > 2s$, $ 0 < \mu < n$.
Subsequently, in \cite{divyadoubly} authors considered the strongly singular (all $ q<0 $ in \eqref{ts}) case to prove the global multiplicity of positive weak solutions using non-smooth analysis. 
Here authors first obtained an elegant comparison theorem for local solutions and used the critical
point theory of non-smooth analysis and the geometry of the energy functional to arrive at the desired results.

In like manner extensive literature is there for Kirchhoff problems due to its vast application in physical and biological models. By considering the effects of change in the length of strings during vibration Kirchhoff in \cite{Kirchhoff} extends the D'Alembert Wave equation. For recent works on existence and multiplicity results for  Kirchhoff type equations we cite the works in \cite{K1, K2, K3, mbr1, mbr2,brw} and references therein. \cite{Kdivya} proved the existence and multiplicity of positive solution of the following critical growth Kirchhoff-Choquard problem using Nehari manifold and Concentration-compactness lemma for the case $1< q< 2$ and for $q = 2$ using Mountain Pass lemma.
\begin{equation*}
	\displaystyle-M(\|\nabla u\|_{L^2}^2)\Delta u = \lambda f(x) |u|^{q-2}u + \left( \int\limits_{\Omega}  \frac{|u(y)|^{2^{*}_{\mu }}}{|x-y|^ \mu}\, dy\right)  |u|^{2^{*}_{\mu }-2}u \;\text{in} \; \Omega,\; u = 0 \; \text{on} \; \partial\Omega,
\end{equation*}
where $\Omega$ is a bounded domain in $\mathbb{R}^N, (N \geq 3)$ with $C^2$ boundary, $M(t) = a + \epsilon^pt^{\theta - 1} , \epsilon > 0$ is small enough,
$0 < \mu < N$, $1 < q \leq 2$, $ a$, $\lambda$, $p$, $\theta$ are positive real numbers such that $p > N-2$, $2_\mu^* = \frac{2N-\mu}{N-2}$, $\theta \in [1, 2_\mu ^{*})$ and f is a continuous real valued sign changing function.
Yu Su and Haibo Chen \cite{chen} investigated about the existence, non-existence and multiplicity of non-trivial solutions of the non-degenerate Kirchhoff-Choquard equations of above type in the whole space.

\noindent The inspiring point of our paper is the work  of Fiscella \cite{fiscella}, where the author studied the degenerate Kirchhoff singular problem. More precisely, they proved multiplicity result using the minimization argument for the first non-trivial solution and approximating the perturbed problem to get second solution for the following Kirchhoff problem with singular and critical non-linearity,
\begin{equation*}
		M_0(\|\nabla u\|_{L^2}^{2})(-\Delta)^s u = \frac{\lambda}{u^\gamma} + u^{2^{*}_{s}-1}, \quad u>0\; \text{in} \; \Omega, \; u = 0 \; \text{on} \;\mathbb{R}^{N} \backslash \Omega,
\end{equation*}
where $M_0(t)=t^{\theta-1},$ $\Omega$ is an open bounded subset with continuous boundary, $N > 2s$, $s, \gamma\in (0, 1)$, $\lambda > 0$, $\theta \in \left( 1, \frac{2_s^{*}}{2}\right)$,  $2_s^{*} = \frac{2N}{N-2s} $ is the fractional critical Sobolev exponent. 

To the best of our knowledge there are no results on Kirchhoff-Choquard term with singularity in both local and non-local case. In this article we study multiplicity result of $(P_\lambda)$ and also comment about the regularity of the problem for the case $\mu < \min\{N, 4s\}$. Using the method of minimization we conclude about the existence of first solution. Moreover, for the second solution we study the perturbed problem obtained by truncating the singular term and then further approximating it to get the second solution for the problem $(P_\lambda)$ using the asymptotic estimates on the minimizers of Hardy-Littlewood-Sobolev non-linearity. With this introduction we state our main results. First we show that any  solution of the problem $(P_\lambda)$:
\begin{Theorem}\label{th1.1}
	Let us suppose $\mu < \min\{N, 4s\}$ and $u\in X_0$ be a solution of the problem $(P_\lambda)$. Then $u\in L^\infty(\Omega)\cap C^{0,s}(\mathbb{R}^N).$
\end{Theorem}
Second result is on the multiplicity of positive solutions:
\begin{Theorem}\label{thm1.1} Let $s \in (0, 1)$, $N > 2s$, $\theta \in [1,2^{*}_{\mu ,s})$, $\gamma\in (0, 1)$ and let $\Omega$ be an open bounded subset of $\mathbb{R}^N$ with $C^2$ boundary. Then there exists $\overline{\lambda} > 0$ such that for any $\lambda\in (0, \overline{\lambda})$ Problem $(P_\lambda)$ has at least two different
	solutions.
\end{Theorem}
\noindent The paper is organized as follows: In section 2, we present some preliminaries on function spaces required for variational settings. In section 3 we see regularity result for the case $ \mu<\min\{N, 4s\}$ and thus prove Theorem \ref{th1.1}. In section 4 we infer the existence of the first positive solution. In section 5 we find a positive solution for the perturbed problem using mountain pass lemma. In section 6 we prove Theorem \ref{thm1.1} and show the existence of the second solution.

	\section{Preliminaries}
	We recall some definitions of function spaces and results that will be required in later sections.
	Consider the functional space $H^s(\mathbb{R}^N)$ as the usual fractional Sobolev space  defined as
\begin{equation*}
	H^s(\mathbb{R}^N) = \left\lbrace u\in L^2(\mathbb{R}^N): \iint\limits_{\mathbb{R}^{2N}} \frac{|u(x)-u(y)|^2}{|x-y|^{N+2s}}\,dxdy < \infty \right\rbrace
\end{equation*}
with the norm 
\begin{equation}\label{eq2.1}
	\|u\|_{H^s(\mathbb{R}^N)} = \|u\|_{L^2(\mathbb{R}^N))} + \displaystyle \left( \iint\limits_{\mathbb{R}^{2N}} \frac{|u(x)-u(y)|^2}{|x-y|^{N+2s}}\,dxdy \right)^{\frac{1}{2}}.
\end{equation}
We define the functional space  $X_0$ as
\begin{equation*}
	X_0 = \left\lbrace u\in H^s(\mathbb{R}^N): u = 0\; a.e.\; in \; \mathbb{R}^N\backslash \Omega  \right\rbrace.
\end{equation*} 
Then, it can be shown that $X_0$ is a Hilbert space with the inner product 
\begin{equation*}
	\langle u, v \rangle = \int\limits_{\mathbb{R}^N}\int\limits_{\mathbb{R}^N}  \dfrac{(u(x)-u(y))(v(x)-v(y))}{|x-y|^{N+2s}}\,dxdy,
\end{equation*}
for $u,v \in X_0$ and thus, the corresponding norm,
\begin{equation*}
	\|u\|_{X_0} = \|u\| = \displaystyle \left( \iint\limits_{\mathbb{R}^{2N}} \frac{|u(x)-u(y)|^2}{|x-y|^{N+2s}}\,dxdy \right)^{\frac{1}{2}}.
\end{equation*} It can be shown that this is equivalent to \eqref{eq2.1} on $X_0$.
\begin{Proposition}\cite{Lieb Analysis}
	\textbf{Hardy–Littlewood–Sobolev inequality}: Let $t$, $r > 1$ and $0 < \mu < N$ with $\frac{1}{t} + \frac{\mu}{N} + \frac{1}{r} = 2$, $f \in L^t(\mathbb{R}^N)$ and $h \in L^r(\mathbb{R}^N)$. Then there exists a sharp constant $C(t, r, \mu, N)$ independent of $f$, $h$ such that
\begin{equation*}
\iint\limits_{\mathbb{R}^{2N}} \dfrac{f(x)h(y)}{|x-y|^ \mu}\,dxdy \leq C(t, r, \mu, N)\|f\|_{L^t(\mathbb{R}^N)}\|h\|_{L^r(\mathbb{R}^N)}.
\end{equation*}

\end{Proposition}

Taking into account that we are looking for positive solutions, we will consider the following problem,
\begin{equation}\label{eq2.2}  \begin{array}{rl}
		\displaystyle \left( \iint\limits_{\mathbb{R}^{2N}} \frac{|u(x)-u(y)|^2}{|x-y|^{N+2s}}\,dxdy \right) ^{\theta - 1}(-\Delta)^s u &= \lambda (u^+)^{-\gamma} + \displaystyle\left( \int\limits_{\Omega} \frac{|u^{+}(y)|^{2^{*}_{\mu ,s}}}{|x-y|^ \mu}\, dy\right)  (u^+)^{2^{*}_{\mu ,s}-1}\;\text{in}\; \Omega\\
		\quad \quad u & = 0\quad \text{in} \; \mathbb{R}^{N}\backslash\Omega.
	\end{array}
\end{equation}
and its weak solution $u \in X_0$ will satisfy,
\begin{equation}\label{eq2.3}
	\displaystyle{\|u\|^{2(\theta-1)}\langle u,\phi \rangle=\lambda \int\limits_{\Omega}(u^+)^{-\gamma}\phi(x) \,dx +\int\limits_{\Omega}\int\limits_{\Omega}  \dfrac{|u^{+}(y)|^{2^{*}_{\mu ,s}}(u^+(x))^{2^{*}_{\mu ,s}-1}\phi(x)}{|x-y|^ \mu}\,dxdy}
\end{equation}
for any $\phi \in X_0$. The energy functional associated with the problem \eqref{eq2.2} is ${J_{\lambda } : X_0 \rightarrow \mathbb{R}}$\; defined as,
\begin{equation*}
	\displaystyle J_\lambda (u)= \dfrac{{\lVert u \rVert}^{2\theta}}{2\theta} - \dfrac{\lambda}{1-\gamma} \int\limits_{\Omega} (u^+)^{1-\gamma} \,dx -\dfrac{1}{2\cdot{2^*_{\mu ,s} }}\int\limits_{\Omega}\int\limits_{\Omega} \dfrac{|u^{+}(x)|^{2^{*}_{\mu ,s}}|u^{+}(y)|^{2^{*}_{\mu ,s}}}{|x-y|^ \mu}\,dxdy.
\end{equation*}
We cannot directly apply critical theory to \eqref{eq2.2}, as the functional $J_\lambda$ is not differentiable on $X_0$ due to the presence of a singular term.
However, with the help of suitable minimization argument we can conclude about the first solution and for the second solution we examine an approximate problem. Now consider the following perturbed problem, for any $n \in \mathbb{N}$:
\begin{equation}\label{pb2.4} (P_{n,\lambda})\;\; \left\{
	\begin{array}{rl}
		  \displaystyle\|u\|^{ 2(\theta -1)}{(-\Delta)}^s u &= \lambda\displaystyle \left(u^+ +\frac{1}{n}\right)^{-\gamma}+\displaystyle\left( \int\limits_{\Omega}\frac{|u^{+}(y)|^{2^*}_{\mu,s}}{|x-y|^\mu}\,dy\right) (u^+)^{2_{\mu,s}^*-1}\quad \text{in} \; \Omega,\\
		u&=0\quad \text{in} \quad\mathbb{R}^N\backslash \Omega,
	\end{array}\right.
\end{equation}
and $u \in X_{0}$ is going to be the weak solution of \eqref{pb2.4} if,
\begin{equation}\label{eq2.5}
	\displaystyle {\|u\|^{2(\theta-1)}\langle u,\phi \rangle=\lambda \int\limits_{\Omega} \left(u^+ +\frac{1}{n}\right)^{-\gamma}\phi(x) \,dx +\int\limits_{\Omega}\int\limits_{\Omega}\frac{|u^{+}(y)|^{2^{*}_{\mu ,s}}(u^+(x))^{2^{*}_{\mu ,s}-1}\phi(x)}{|x-y|^ \mu}\,dxdy},
\end{equation} for all $\phi$ in $X_0$.
Now, the critical point of $J_{n,\lambda}:X_0 \rightarrow \mathbb{R}$, are the weak solutions of \eqref{pb2.4}, where $J_{n,\lambda}$ defined as,
\begin{equation}\label{eq2.6}
		\begin{aligned}
		\displaystyle J_{n,\lambda}(u) = & \frac{\|u\|^ {2\theta}}{2\theta}-\frac{\lambda}{1-\gamma}\int\limits_{\Omega}\left[ {\left(u^+ +\frac{1}{n}\right)}^{1-\gamma}-\left(\frac{1}{n}\right)^{1-\gamma}\right] \,dx\\ &-\frac{1}{2\cdot{2^*}_{\mu,s}}\int\limits_{\Omega}\int\limits_{\Omega}\frac{|u^{+}(y)|^{2^*_{\mu,s}}|u^{+}(x)|^{2_{\mu,s}^*}}{|x-y|^\mu}\,dxdy.
	\end{aligned}
\end{equation}
From the embedding results, we conclude that $X_0$ is continuously embedded in $L^{p}(\Omega)$ when $1 \leq p \leq 2_s^{*}$. Also the embedding is compact for $1 \leq p < 2_s^{*}$, but not for the case $p = 2_s^{*}$. We define
 the best constant for the embedding $X_0$ into $L^{2_s^{*}}(\mathbb{R}^N)$ as, 
 \begin{equation}\label{eq2.7}
 	S_s =\inf\limits_{u\in X_0\backslash \{0\}}\left\lbrace\iint\limits_{\mathbb{R}^{2N}} \frac{|u(x)-u(y)|^2}{|x-y|^{N+2s}}\,dxdy:  \int\limits_{\mathbb{R^N}}|u|^{2_s^{*}} = 1 \right\rbrace.
 \end{equation}
 Consequently, we define
 \begin{equation}\label{eq2.8}
 	S_s^H = \inf\limits_{u\in X_0\backslash \{0\}}\left\lbrace\iint\limits_{\mathbb{R}^{2N}} \frac{|u(x)-u(y)|^2}{|x-y|^{N+2s}}\,dxdy:  \iint\limits_{\mathbb{R}^{2N}} \dfrac{|u(x)|^{2^{*}_{\mu ,s}}|u(y)|^{2^{*}_{\mu ,s}}}{|x-y|^ \mu}\,dxdy = 1\right\rbrace. 
 \end{equation}
\begin{Lemma} \cite{Tuhina}
	The constant $S_s^H$ is achieved by u if and only if u is of the form\\
	$C\left( \frac{t}{t^2 + |x-x_0|^2}\right) ^\frac{N-2s}{2}$, $x \in \mathbb{R}^N$, 
	for some $x_0 \in \mathbb{R}^N, C \text{and}\; t > 0.$ Moreover,
	$S_s^H = \frac{S_s}{{C(N, \mu)}^\frac{1}{2_{\mu, s}^{*}}}$.
\end{Lemma}

Consider the family of functions ${U_\epsilon}$, where $U_\epsilon$ is defined as
\begin{equation*}
	U_\epsilon = \epsilon^{-\frac{N-2s}{2}}u^{*}\left( \frac{x}{\epsilon}\right), \; x\in \mathbb{R}^{N}, \epsilon > 0,
\end{equation*}
\begin{equation*}
	u^{*}(x) = \overline{u}\left( \frac{x}{{S_s}^\frac{1}{2s}}\right) , \; \overline{u}(x) = \frac{\tilde{u}(x)}{\|\tilde u\|_{L^{2_s^{*}}(\mathbb{R}^N)}} \; \text{and}\; \tilde{u}(x) = \alpha(\beta^2 + |x|^2)^{-\frac{N-2s}{2}},
	\end{equation*}
with $\alpha \in \mathbb{R}\backslash\{0\}$ and $\beta > 0$ are fixed constants.
Then for each $\epsilon > 0, \; U_\epsilon $ satisfies
\begin{equation*}
	(-\Delta)^s u = |u|^{2_s^{*}-2}u \quad in \; \mathbb{R}^N,
\end{equation*}
and the equality,
\begin{equation*}
	\iint\limits_{\mathbb{R}^{2N}} \frac{|U_\epsilon(x)-U_\epsilon(y)|^2}{|x-y|^{N+2s}}\,dxdy = \int\limits_{\mathbb{R^N}}|U_\epsilon|^{2_s^{*}} = {S_s}^\frac{N}{2s}. 
\end{equation*}

Without loss of generality, we assume $0 \in \Omega$  and fix $\delta >0$ such that $B_{4\delta} \subset \Omega$. Let   a $\eta \in C^{\infty}(\mathbb{R}^N)$ be a cut off function  such that
\begin{equation*}
\eta =	\begin{cases}
		 1 & \quad B_{\delta},\\
		0 & \quad \mathbb{R}^N\backslash B_{2\delta},
	\end{cases}
\end{equation*}
and for each $\epsilon > 0$, let  $u_\epsilon$ be defined as \begin{equation}\label{eq2.9}
u_\epsilon(x) = \eta(x)U_\epsilon(x)\quad for \; x\in \mathbb{R}^N.
\end{equation}
Then we have the following proposition from \cite[Propositions 2.7 and 2.8]{Tuhina}.
\begin{Proposition}\label{Prop2.3}
	Let $s>0 $ and $N > 2s$. Then, the following estimates hold true as $\epsilon \to 0$
	\begin{enumerate} 
 	\item $\displaystyle\iint\limits_{\mathbb{R}^{2N}} \frac{|u_\epsilon(x)-u_\epsilon(y)|^2}{|x-y|^{N+2s}}\,dxdy  \leq  S_s^{\frac{N}{2s}} + O(\epsilon^{N-2s}) = C^\frac{N(N-2s)}{2s(2N-\mu)}(S^H_s)^\frac{N}{2s}+ O(\epsilon^{N-2s})$,
    \item $\displaystyle \int\limits_{\Omega}\int\limits_{\Omega} \frac{|u_\epsilon(x)|^{2^{*}_{\mu ,s}}|u_\epsilon(y)|^{2^{*}_{\mu ,s}}}{|x-y|^ \mu}\,dxdy \geq C^\frac{N}{2s}(S^H_s)^\frac{2N-\mu}{2s}- O(\epsilon^{N})$. 
\end{enumerate}
\end{Proposition}
Throughout the paper we will use the notation,
\begin{equation*}
	\|u\|_{N L}^{2\cdot2^{*}_{\mu ,s}} := \int\limits_{\Omega}\int\limits_{\Omega} \frac{|u(x)|^{2^{*}_{\mu ,s}}|u(y)|^{2^{*}_{\mu ,s}}}{|x-y|^ \mu}\,dxdy.
\end{equation*}
We recall that ${(u_k)}_k \subset X_0$ is a Palais–Smale sequence for $C^1$ functional $J : X_0 \rightarrow \mathbb{R}$ at level $c \in \mathbb{R}$ if,
\begin{equation}\label{eq4.2}
	J(u_k) \rightarrow c \quad \text{and} \quad J'(u_k) \rightarrow 0 \quad \text{in} \; X'_0 \quad \text{as} \; k \to  \infty.
\end{equation}
We say that $J$ satisfies the Palais–Smale condition at level $c$ if any Palais–Smale sequence ${(u_k)}_k$ at level $c$
admits a convergent subsequence in $X_0$.

\section{Regularity of solutions of $(P_\lambda)$}
Let us define $d :\overline{\Omega} \rightarrow \mathbb{R}^+ $ by $d(x) := dist(x,\mathbb{R}^N\backslash\Omega)$. Now given any $\phi \in C^0(\overline{\Omega})$ such that $\phi > 0$ in $\Omega$ we define
\begin{equation*}
C_\phi(\Omega) := \left\lbrace u \in C_0(\overline{\Omega}) :\; \text{exists}\; c \geq 0\;\text{such that} \;|u(x)| \leq c\phi(x),\; \text{for all}\; x\in \Omega\right\rbrace,
\end{equation*}
with the usual norm $\displaystyle{\|\frac{u}{\phi}\|}_{{L_\infty}(\Omega)}$. We define
the following open convex subset of $C_\phi(\Omega)$ as
\begin{equation*}
C_\phi^+(\Omega) := \left\lbrace u \in C_\phi(\Omega) :\; \inf\limits_{x \in \Omega}\frac{u(x)}{\phi(x)} > 0\right\rbrace.
\end{equation*}
Let $\phi_1$ be the first positive normalized eigen-function of $(-\Delta)^s$ in $X_0$. From \cite[Proposition 1.1, Theorem 1.2]{ros-oton } we recall that $\phi_1\in C^{0,s}(\mathbb{R}^N)\cap C_{d^s}^+(\Omega).$
Next we will show that the solutions of $(P_\lambda)$ are bounded and H\"older continuous in the case $\mu < 4s$. 

{\bf Proof of Theorem \ref{th1.1}:} Let $\zeta\in C^1(\mathbb{R},[0,1])$ be a cut-off  function 
 such that  $ \zeta(t)=
0,$ when $ t\leq0$ and equal to $1$, for  $t\geq1.$
Then 
$	\zeta'(t)\geq 0,\; 0\leq t\leq 1. $
Now, for $\epsilon>0$, define
\begin{center}$\zeta_\epsilon(t)= \zeta(\dfrac{t-1}{\epsilon}).$
\end{center} 
\begin{equation*}
\zeta'_\epsilon(t)\geq 0, \quad \text{for all} \; 1\leq t \leq 1+\epsilon.
\end{equation*}
Let $u$ be a solution of problem $(P_\lambda)$ and so positive, $\kappa\in [C_c^\infty(\Omega)]^+ ;\ \eta = (\zeta_\epsilon \circ u)\kappa \in C_c^\infty$ is an appropriate test function. We get,
\begin{equation*}
\begin{aligned}
\|u\|^{2(\theta-1)}\int\limits_{\mathbb{R}^N}\int\limits_{\mathbb{R}^N}  \dfrac{(u(x)-u(y))(\eta(x)-\eta(y))}{|x-y|^{N+2s}} = & \lambda\int\limits_{\Omega} u^{-\gamma}(x)\eta(x)\,dx\\ &+ \int\limits_{\Omega}\int\limits_{\Omega}\dfrac{(u(y))^{{2}_{\mu,s}^{*}}(u(x))^{{2}_{\mu,s}^{*}-1}\eta(x)}{|x-y|^\mu}dxdy 
\end{aligned}
\end{equation*}
as $\epsilon\to 0$, take $ v= (u-1)^+$,
\begin{equation}\label{eq6.1}
\begin{aligned}
\|u\|^{2(\theta-1)}\int\limits_{\mathbb{R}^N}\int\limits_{\mathbb{R}^N}  \dfrac{(v(x)-v(y))(\kappa(x)-\kappa(y))}{|x-y|^{N+2s}}\,dxdy \leq & \lambda\int\limits_{\Omega} \kappa(x)\,dx \\& + \int\limits_{\Omega}\int\limits_{\Omega}\dfrac{(v(y))^{{2}_{\mu,s}^{*}}(v(x))^{{2}_{\mu,s}^{*}-1}\kappa(x)}{|x-y|^\mu}\,dxdy \\&+ \int\limits_{\Omega}\int\limits_{\Omega}\dfrac{1\cdot(v(x))^{{2}_{\mu,s}^{*}-1}\kappa(x)}{|x-y|^\mu}\,dxdy.
\end{aligned}
\end{equation}
Let us define, $\beta> 1$ and $T > 0$
\begin{center}
	$\phi(t) = \begin{cases}
	0, & t\leq 0,\\ t^\beta, & 0 < t \leq T, \\ \beta T^{\beta-1}t -(\beta-1)T^\beta, & t > T.
	\end{cases}$
\end{center}
Then,
\begin{center}
	$\phi ^{'}(t) = \begin{cases}
	0, & t\leq 0,\\ \beta t^{\beta-1}, & 0 < t \leq T, \\ \beta T^{\beta-1}, & t > T.
	\end{cases}$
\end{center}
As $v\geq 0$, it's easier to see that $\phi(v)\geq 0$, $\phi'(v)\geq 0$ and  $v\phi'(v) \leq \beta\phi(v)$.
Since $\phi$ is a Lipschitz function and $v \in H_0^s(\Omega)$. This implies $\phi(v) \in H_0^s(\Omega)$ and
\begin{equation*}
\|\phi(v)\|_{X_0}^2=\int\limits_{\mathbb{R}^N}\int\limits_{\mathbb{R}^N}   \dfrac{{|\phi(v(x))-\phi(v(y))|}^{2}}{|x-y|^{N+2s}}\,dxdy   \leq K \|v\|_{X_0}^2.
\end{equation*}
By Sobolev embedding we have,
\begin{equation}\label{eq6.2}
S_s \|\phi(v)\|^{2}_{L^{2_s^*}} \leq \|\phi(v)\|^2		
 = \int\limits_{\mathbb{R}^N} \phi(v) (-\Delta)^{s} \phi (v).
\end{equation}
Since $\phi $ is a convex function, we also have
\begin{equation*}
(-\Delta)^{s} \phi (v) \leq \phi ^{'}(v)(-\Delta)^{s}(v)\quad \text{weakly in} \;\Omega,
\end{equation*}
from \eqref{eq6.1} we see\begin{align*}
\|u\|^{2(\theta-1)}(-\Delta)^s\phi(v) &\leq \|u\|^{2(\theta-1)}\phi ^{'}(v)(-\Delta)^{s}(v)\\
&\leq \phi ^{'}(v) \left[ \lambda + \left( \int\limits_{\Omega}\frac{|v(y)|^{{2}_{\mu,s}^{*}}}{|x-y|^\mu}\,dy\right)  (v)^{{2}_{\mu,s}^{*}-1} + \left( \int\limits_{\Omega}\frac{dy}{|x-y|^\mu}\right) (v)^{{2}_{\mu,s}^{*}-1}\right],
\end{align*}
as $\phi(v)\phi ^{'}(v) \in X_0$ and using \eqref{eq6.2} we get,
\begin{align*}
\|u\|^{2(\theta-1)} S_s \|\phi(v)\|^{2}_{L^{2_s^*}} \leq& \|u\|^{2(\theta-1)} \int\limits_{\Omega} \phi(v) (-\Delta)^{s} \phi (v)\\
 \leq& \lambda\int\limits_{\Omega} \phi(v)\phi ^{'}(v)\,dx+ \int\limits_{\Omega}\int\limits_{\Omega}\dfrac{(v(y))^{{2}_{\mu,s}^{*}}(v(x))^{{2}_{\mu,s}^{*}-1}\phi(v)\phi ^{'}(v)}{|x-y|^\mu}\,dxdy\\& + \int\limits_{\Omega}\int\limits_{\Omega}\dfrac{1\cdot(v(x))^{{2}_{\mu,s}^{*}-1}\phi(v)\phi ^{'}(v)}{|x-y|^\mu}\,dxdy,
\end{align*}
Moreover, by the relation, $v \phi^{'}(v) \leq \beta \phi(v)$ and $\phi ^{'}(v) \leq \beta [1+ \phi(v)]$
\begin{equation*}
\begin{aligned}
	\|\phi(v)\|^{2}_{L^{2_s^*}} \leq& \frac{\beta}{\|u\|^{2(\theta-1)} S_s} \left[ \lambda\int\limits_{\Omega} [\phi(v) + \phi^{2}(v)] \, dx + \int\limits_{\Omega}\int\limits_{\Omega}\frac{(v(y))^{{2}_{\mu,s}^{*}}(v(x))^{{2}_{\mu,s}^{*}-2}\phi^{2}(v)}{|x-y|^\mu}\,dxdy\right]\\  &+\frac{\beta}{\|u\|^{2(\theta-1)} S_s} \left[ \int\limits_{\Omega}\int\limits_{\Omega}\frac{1\cdot(v(x))^{{2}_{\mu,s}^{*}-2}\phi^{2}(v)}{|x-y|^\mu}\,dxdy  \right].
\end{aligned}
\end{equation*}
Now following the proof of Theorem 1.1 of \cite{DS} we obtain  $u\in L^{\infty}(\Omega)$.

To show that $u \in C^{0,s}(\mathbb{R}^N)$, we note that $u$ satisfies the equation 
\begin{equation*}
	\displaystyle(-\Delta)^s u =\frac{1}{\|u\|^{2(\theta-1)}} \left( \frac{\lambda}{u^\gamma} + \left( \int\limits_{\Omega} \frac{|u(y)|^{2^{*}_{\mu ,s}}}{|x-y|^ \mu}\, dy\right)  |u|^{2^{*}_{\mu ,s}-1}\right) := g.\, (say)
\end{equation*}
Then as $u \in X_0\cap L^{\infty}(\Omega)$, we get  $g \in L_{loc}^{\infty}$.  thus from \cite{brasco} we have \begin{equation}\label{eq2.21}
u \in C_{loc}^{0,\alpha}\quad \text{where}\; \alpha < \min\{1, 2s\}.
\end{equation}
Further from \cite[Theorem 1.2]{adimurthi} (for the case $\beta = 0$) let $\underline{u}$ be the unique solution of \begin{equation*}
	\|u\|^{2(\theta-1)}(-\Delta^s v) = v^{-\gamma},
\end{equation*}
and by \cite[Remark 1.5]{adimurthi} ( for $\beta = 0$ and $c =\displaystyle \left|\left(\int \limits_{\Omega}\frac{|u(y)|^{{2}_{\mu,s}^{*}}}{|x-y|^{\mu}}\,dy\right) (u)^{{2}_{\mu,s}^{*}-1} \right|_{L^{\infty}}) $, let $\overline{u}$  be the unique solution of,
\begin{equation*}
	\|u\|^{2(\theta-1)}(-\Delta^s v) = v^{-\gamma} +c,
\end{equation*}
Then $\underline{u}$ and $\overline{u}$ forms the sub-solution and the super-solution of problem$(P_\lambda)$ respectively. We claim that $\underline{u} \leq u \leq \overline{u}$\, a.e. in $\Omega$. First we prove $\underline{u} \leq u$ a.e. in $\Omega$, which follows from \cite[Proposition 4.5]{divyareg}, that is
we assume by contradiction that $\underline{u} > u$. Then for any $u \in X_0$ we have,
 \begin{equation*}
	(u(x)- u(y))(u^+(x)-u^+(y))\geq {|u^+(x)-u^+(y)|}^2
\end{equation*}
This implies, \begin{equation}\label{eq2.22}
	\|u^+\|^2 \leq \int\limits_{\mathbb{R}^N}\int\limits_{\mathbb{R}^N}  \dfrac{(u(x)-u(y))(u^+(x)-u^+(y))}{|x-y|^{N+2s}}\,dxdy
\end{equation}
Testing $\|u\|^{2(\theta-1)}\left( (-\Delta)^s\underline{u} - (-\Delta)^s u \right) \leq {\underline{u}}^{-\gamma} - {u}^{-\gamma} $ with ${(\underline{u}-u)}^+$ we deduce,
\begin{equation*}
\|u\|^{2(\theta-1)}\langle \underline{u}-u, (\underline{u}-u)^+\rangle \leq \int\limits_{\Omega}\left(  {\underline{u}}^{-\gamma} - {u}^{-\gamma}\right) \left( \underline{u}-u\right)^+\,dx\, \leq 0,
\end{equation*}
by \eqref{eq2.22} and using the fact that $u$ is a positive weak solution we get,
\begin{equation*}
	0\leq \|u\|^{2(\theta-1)}\|\left( \underline{u}-u\right)^+\|^2 \leq \int\limits_{\Omega}\left(  {\underline{u}}^{-\gamma} - {u}^{-\gamma}\right) \left( \underline{u}-u\right)^+\,dx\, \leq 0. 
\end{equation*}
Thus, a contradiction. Hence, $\underline{u} \leq u$.
In the similar lines one can prove $\overline{u} \geq u$. 
 Moreover, from \cite[Theorem 1.2, Remark 1.5]{adimurthi}\, $\underline{u},\, \overline{u} \in C_{\phi_1}^+$, where $\phi_1 \in C_{d^s}^+$. This implies, there exist positive constants say $C_1$ and $C_2$ such that
\begin{equation}\label{eq2.23}
C_1 d^s \leq \underline{u} \leq u \leq \overline{u} \leq C_2 d^s,
\end{equation}
So by \eqref{eq2.21} and \eqref{eq2.23} we have $u\in C^0(\overline{\Omega})$, hence it is a classical solution (see \cite[Definition 2]{adimurthi}). Now by following the arguments as in Theorem 1.4 of \cite{adimurthi}  we get
\[ \frac{|u(x)-u(y)|}{|x-y|^s} \le \frac{|u(x)}{|x-y|^s}+\frac{|u(y)|}{|x-y|^s}\leq  2^s \left( \frac{u(x)}{d(x)^s}+\frac{u(y)}{d^s(y)}\right)\] 
for all $x,y \in \Omega_\eta=\{ x\in \Omega, d(x)<\eta\}$ with $|x-y|\ge \max\left\{\frac{d(x)}{2},\frac{d(y)}{2} \right \}.$ Combining with interior regularity, we get
$u \in C^{0,s}(\mathbb{R}^N)$.\QED
   
	\section{Existence of first solution}
	In this section we show the existence of solution for \eqref{eq2.2} by showing the existence of minimizer for $J_\lambda$ over small ball around $0$. 
	\begin{Lemma}\label{Lemma1}
	There exist numbers $\displaystyle \rho  \in (0, 1]$, $\lambda_0 = \lambda_0(\rho) >0$ and $\alpha = \alpha(\rho) >0 $ such that $J_\lambda (u)\geq  \alpha $ for any $u \in X_0$ ,
	with  $\|u\|  =\rho$ , and for any $\lambda \in (0, \lambda_0]$. Furthermore, set\\
	$ \displaystyle m_\lambda =\inf \{J_\lambda (u) : u \in \overline{\rm B_\rho}\}$,
		where $\overline{\rm B_\rho} = \{ u \in X_0 : \lVert u \rVert \leq \rho \}$. Then ${m_\lambda < 0}$ for any ${\lambda \in (0, \lambda_0]}$.
	\end{Lemma}
\proof  From Hardy-Littlewod-Sobolev inequality, we have

\begin{equation}\label{eq3.2}
		\int\limits_{\Omega}\int\limits_{\Omega} \dfrac{|u^{+}(x)|^{2^{*}_{\mu ,s}}|u^{+}(y)|^{2^{*}_{\mu ,s}}}{|x-y|^ \mu} \,dxdy \leq C(\mu,N)S_s^{-2^{*}_{\mu ,s}} {\lVert u\rVert}^{2\cdot{2^*_{\mu ,s} }}
	\end{equation}
	Then using the H\"older inequality and \eqref{eq3.2} we get,
		\begin{equation*} J_\lambda (u)\geq \dfrac{{\lVert u \rVert}^{2\theta}}{2\theta}-\frac{\lambda}{1-\gamma}|\Omega|^{\frac{{2^*_s}-1+\gamma}{2^*_s}} S_s^{\frac{-(1-\gamma)}{2}}{\lVert u\rVert}^{1-\gamma}-\frac{ C(\mu,N)S_s^{-2^{*}_{\mu ,s}} {\lVert u\rVert}^{2\cdot{2^*_{\mu ,s} }}}{2\cdot{2^*_{\mu ,s} }}.
	\end{equation*}
	 
	Since $1-\gamma < 1< 2\theta <2\cdot{2_{\mu ,s}^{*}}$, the function 
	\begin{equation*}
		\eta(t) = {\dfrac{t^{2\theta-1+\gamma}}{2\theta}}-{\dfrac{ C(\mu,N)S_s^{-2^{*}_{\mu ,s}} }{2\cdot{2^*_{\mu ,s} }}}{t}^{2\cdot{2^*_{\mu ,s} }-1+\gamma}, \quad t \in [0,1],
	\end{equation*}
	attains maxima at some $\rho \in(0,1]$ small enough, that is, ${\max\limits_{t\in [0,1]} \eta(t) =\eta(\rho) > 0}$. Letting
	$$\lambda_0 =\frac{(1-\gamma)S_s^{\frac{(1-\gamma)}{2}}}{2|\Omega|^{\frac{{2^*_s}-1+\gamma}{2^*_s}}}\eta(\rho),$$
	it follows that there exists a number $\alpha > 0$ such that for any $u\in X_0$ with $\lVert u\rVert = \rho \leq 1$ and for any $\lambda \leq \lambda_0$, we get
	\begin{equation*}
		\displaystyle{J_\lambda (u)\geq {\rho^{1-\gamma} } {\eta(\rho)}/2=\alpha >0}.
	\end{equation*}
	 Since $1-\gamma < 1< 2\theta <2\cdot{2_{\mu ,s}^{*}}$, so for a fixed $v\in X_0$ with $v^+\neq 0$ and $t$ small enough we have,
	\begin{equation*}
		J_\lambda(tv) = \dfrac{t^{2\theta}}{2\theta}{\lVert v \rVert}^{2\theta} - \dfrac{\lambda t^{1-\gamma}}{1-\gamma} \int\limits_{\Omega} (v^+)^{1-\gamma} \,dx -\dfrac{t^{2\cdot{2^*_{\mu ,s}}}}{2\cdot{2^*_{\mu ,s} }}\int\limits_{\Omega}\int\limits_{\Omega} \dfrac{|v^{+}(x)|^{2^{*}_{\mu ,s}}|v^{+}(y)|^{2^{*}_{\mu ,s}}}{|x-y|^ \mu}\,dxdy < 0,
	\end{equation*}
	 this implies for $\|u\|$ small enough, $m_\lambda < 0$ for $\lambda \in (0, \lambda_0]$.
	\QED
   
	\begin{Theorem}\label{thm3.2} Let $\lambda_0$ be given as in Lemma \ref{Lemma1}. Then for any $\lambda \in (0, \lambda_0]$ problem $(P_\lambda)$ has a solution $u_0 \in X_0$
	with $J_\lambda (u_0) < 0$.
    \end{Theorem}

	\proof First we claim that there exists $u_0\in \overline{\rm B_\rho } $ such that $J_\lambda(u_0) = m_\lambda < 0$ where $\lambda \in (0, \lambda_0]$ and let $\rho$ be as given in Lemma \ref{Lemma1}. Definition of $m_\lambda$, implies that there exists a minimizing sequence say ${(u_k)}_k \subset \overline{\rm B_\rho}$ such that
   \begin{equation}\label{eq3.3}
	\lim \limits_{k \to \infty} J_\lambda(u_k) = m_\lambda.
   \end{equation}
	Clearly ${(u_k)}_k$ is a bounded sequence in $X_0$, up to a sub-sequence,  there exists a function $u_0 \in X_0$ such that, as $k \to \infty $, we have
$u_k \rightharpoonup u_0 $ weakly in $X_0, u_k \rightharpoonup u_0$ weakly in $ L^{2_s^*}(\Omega), u_k \rightarrow u_0$ strongly in $ L^p(\Omega)$ for all $ p\in[1,2_s^*), u_k \rightarrow u_0$ a.e in $\Omega.$\, Since $\gamma \in (0,1)$, for any $k \in \mathbb{N}$ we have
	\begin{equation}\label{eq3.5}
			\lim \limits_{k \to \infty} \int\limits_{\Omega}(u_k^+)^{1-\gamma} \,dx =\int\limits_{\Omega}(u_0^+)^{1-\gamma} \,dx.
		\end{equation}
  We set ${w_k = u_k-u_0}$; by \cite[Theorem 2]{7thm2} and \cite[Lemma 2.2]{M.yang} one has,
	\begin{equation}\label{eq3.6}
		\|u_k\|^2= \|w_k\|^2+\|u_0\|^2+ o(1),\; \;
			\|u_k\|_{N L}^{2\cdot2^{*}_{\mu ,s}} =\|w_k\|_{N L}^{2\cdot2^{*}_{\mu ,s}} + \|u_0\|_{N L}^{2\cdot2^{*}_{\mu ,s}} +o(1).
	    \end{equation}
		By \eqref{eq3.6}, we have $w_k \in \overline{\rm B_\rho}$ for $k$ sufficiently large 
		and from this, since $\rho \leq 1$ we have
			\begin{equation}\label{eq3.8}
				\dfrac{{\lVert w_k \rVert}^{2\theta}}{2\theta} -\dfrac{\|w_k^{+}\|_{N L}^{2\cdot2^{*}_{\mu ,s}}}{2\cdot{2^*_{\mu ,s} }} > 0.
			\end{equation}
		Thus, by \eqref{eq3.3}, \eqref{eq3.5}, \eqref{eq3.8} it follows that, as $k \to \infty$
	    \begin{align*}
	    	m_\lambda &= J_\lambda (u_k) + o(1)\\
	    	&=\dfrac{({{\lVert w_k \rVert}^{2}+{\lVert u_0 \rVert}^{2}})^{\theta}}{2\theta}-\dfrac{\lambda }{1-\gamma} \int\limits_{\Omega} (u_0^+)^{1-\gamma} \,dx-\dfrac{\|w_k^{+}\|_{N L}^{2\cdot2^{*}_{\mu ,s}}}{2\cdot{2^*_{\mu ,s} }}-\dfrac{\|u_0^{+}\|_{N L}^{2\cdot2^{*}_{\mu ,s}}}{2\cdot{2^*_{\mu ,s} }} + o(1)\\
	    	&\geq J_\lambda (u_0)+\dfrac{{\lVert w_k \rVert}^{2\theta}}{2\theta}-\dfrac{\|w_k^{+}\|_{N L}^{2\cdot2^{*}_{\mu ,s}}}{2\cdot{2^*_{\mu ,s} }} +o(1)\\
	    	&\geq J_\lambda (u_0) + o(1)
	    	\geq m_\lambda.
	    \end{align*}
    Noting that $\overline{\rm B_\rho }$ is closed convex set, thus $u_0\in \overline{\rm B_\rho } $. Hence, $u_0$ is a local minimizer for $J_\lambda$, with $J_\lambda(u_0) = m_\lambda < 0$, which implies that $u_0$ is nontrivial.\\
	Next, it remains to prove that $u_0$ is a positive solution of \eqref{eq2.2}. As $u_0$ is a local minimizer for $J_\lambda$, thus for any $\psi\in X_0$, with $\psi \geq 0$ a.e. in $\mathbb{R^N}$ and $t > 0$ be small enough so that $u_0 + t\psi \in \overline{\rm B_\rho}$\; we have
	\begin{align*}
		0\leq & J_\lambda(u_0 + t\psi)- J_\lambda (u_0)\\
		=& \displaystyle \frac{{\lVert u_0 + t\psi \rVert}^{2\theta}-{\lVert u_0 \rVert}^{2\theta}}{2\theta}-\dfrac{\lambda }{1-\gamma} \int\limits_{\Omega}\left[ ({(u_0 + t\psi)}^+)^{1-\gamma}-(u_0^+)^{1-\gamma}\right]  \,dx \\&-\dfrac{1}{2\cdot{2^*_{\mu ,s}}} \left( \|(u_0 + t\psi)^{+}\|_{N L}^{2\cdot2^{*}_{\mu ,s}}-\|u_0^{+}\|_{N L}^{2\cdot2^{*}_{\mu ,s}}\right) .
	\end{align*}
    Note that,
    \begin{equation*}
    	\dfrac{1 }{1-\gamma} \frac{({(u_0 + t\psi)}^+)^{1-\gamma}- (u_0^+)^{1-\gamma}}{t} = ({(u_0 + \xi t\psi)}^+)^{-\gamma}\psi,
    \end{equation*}
    a.e. in $\Omega$ with $ \xi \in (0, 1)$ and $({(u_0 + \xi t\psi)}^+)^{-\gamma}\psi \to ({(u_0)}^+)^{-\gamma}\psi$\;  a.e. in $\Omega$ as $t\to 0^+$. From the Fatou lemma, it follows that
    \begin{equation}\label{eq3.9}
    	\lambda \int\limits_{\Omega} (u_0^+)^{-\gamma}\psi \,dx \leq \liminf \limits_{t \to 0+} \dfrac{\lambda }{1-\gamma} \int\limits_{\Omega} \frac{({(u_0 + t\psi)}^+)^{1-\gamma}- (u_0^+)^{1-\gamma}}{t} \,dx.
    \end{equation}
    So dividing by $t > 0$ and taking $t \to 0^+$, we get
	\begin{equation}\label{eq3.10}
		\begin{aligned}
			\liminf \limits_{t \to 0+} \dfrac{\lambda }{1-\gamma} \int\limits_{\Omega} \frac{({(u_0 + t\psi)}^+)^{1-\gamma}- (u_0^+)^{1-\gamma}}{t}dx \leq& \|u_0\|^{2(\theta -1)}\langle u_0,\phi \rangle\\&-\int\limits_{\Omega}\int\limits_{\Omega} \dfrac{|u_0^{+}(y)|^{2^{*}_{\mu ,s}}|u_0^{+}(x)|^{2^{*}_{\mu ,s}-2}u_0^{+}(x)\psi (x)}{|x-y|^ \mu}.
		\end{aligned}
		\end{equation}

	Therefore, by \eqref{eq3.9}, \eqref{eq3.10} for any $\psi \in X_0$ with $\psi \geq 0$ a.e. in $\mathbb{R}^{N}$,
		\begin{equation}\label{eq3.11}
			0 \leq \|u_0\|^{2(\theta -1)}\langle u_0,\psi \rangle-\lambda \int\limits_{\Omega} (u_0^+)^{-\gamma}\psi \,dx -\int\limits_{\Omega}\int\limits_{\Omega}\dfrac{|u_0^{+}(y)|^{2^{*}_{\mu ,s}}(u_0^+(x))^{2^{*}_{\mu ,s}-1}\psi (x)}{|x-y|^ \mu}\,dxdy.
		\end{equation}
			Since $J_\lambda (u_0) < 0$ and by Lemma \ref{Lemma1}, we have $u_0 \in \overline{\rm B_\rho }$. Hence, there exists $ \delta \in (0, 1)$ such that $(1+t)u_0 \in \overline{\rm B_\rho }$
			for any $t \in [-\delta, \delta]$. Define $I_\lambda (t) = J_\lambda((1 + t)u_0)$. Clearly, the functional $I_\lambda$ has a minimum at $t = 0$ i.e.
			 \begin{equation}\label{eq3.12}
			 	I_\lambda^{'}(0) = \|u_0\|^{2\theta}-\lambda \int\limits_{\Omega} (u_0^+)^{1-\gamma} \,dx -\int\limits_{\Omega}\int\limits_{\Omega}\dfrac{|u_0^{+}(y)|^{2^{*}_{\mu ,s}}|u_0^+(x)|^{2^{*}_{\mu ,s}}}{|x-y|^ \mu}\,dxdy = 0.
			 \end{equation}
		Suppose for any $\epsilon > 0$ and $\phi \in X_0$, we define $\psi _\epsilon = u_0^+ +\epsilon \phi$. By \eqref{eq3.11} we have
	    \begin{equation}\label{eq3.13}
		\begin{aligned}[t]
      0 \leq& \|u_0\|^{2(\theta -1)}\langle u_0,\psi_\epsilon^+ \rangle-\lambda \int\limits_{\Omega}(u_0^+)^{-\gamma}\psi_\epsilon^+ \,dx -\int\limits_{\Omega}\int\limits_{\Omega}\dfrac{|u_0^{+}(y)|^{2^{*}_{\mu ,s}}(u_0^+(x))^{2^{*}_{\mu ,s}-1}\psi_\epsilon^+ (x)}{|x-y|^ \mu}\,dxdy\\
      =&\|u_0\|^{2(\theta -1)}\langle u_0,\psi_\epsilon +\psi_\epsilon^- \rangle-\lambda \int \limits_{\Omega}(u_0^+)^{-\gamma}(\psi_\epsilon +\psi_\epsilon^-) \,dx\\ &-\int\limits_{\Omega}\int\limits_{\Omega}\dfrac{|u_0^{+}(y)|^{2^{*}_{\mu ,s}}(u_0^+(x))^{2^{*}_{\mu ,s}-1}(\psi_\epsilon +\psi_\epsilon^-) (x)}{|x-y|^ \mu}\,dxdy.
  \end{aligned}
		\end{equation}
    Also note that for $x, y \in\mathbb{R^N}$
		\begin{equation}\label{eq3.14}
			\begin{aligned}
(u_0(x)- u_0(y))(u_0^-(x)-u_0^-(y)) &= -u_0^+(x)u_0^-(y)-u_0^-(x)u_0^+(y)-{[u_0^-(x)-u_0^-(y)]}^2\\
&\leq -{|u_0^-(x)-u_0^-(y)|}^2,
			\end{aligned}
		\end{equation}
		and from the previous inequality we also get,
		$$(u_0(x)- u_0(y))(u_0^+(x)-u_0^+(y))\leq {|u_0(x)-u_0(y)|}^2.$$
	This implies,
			\begin{equation}\label{eq3.15}
			\begin{aligned}[t]
			\langle u_0,\psi_\epsilon +\psi_\epsilon^- \rangle =& \iint\limits_{\mathbb{R}^{2N}} \dfrac{(u_0(x)-u_0(y))(\psi_\epsilon(x) +\psi_\epsilon^-(x)-\psi_\epsilon(y) -\psi_\epsilon^-(y))}{|x-y|^{N+2s}}\,dxdy\\
			=&\iint\limits_{\mathbb{R}^{2N}}\dfrac{|u_0(x)-u_0(y)|^2}{|x-y|^ {N+2s}}\,dxdy+\epsilon\iint\limits_{\mathbb{R}^{2N}} \dfrac{(u_0(x)-u_0(y))(\phi(x)-\phi(y))}{|x-y|^ {N+2s}}\,dxdy\\ &+\iint\limits_{\mathbb{R}^{2N}}\dfrac{(u_0(x)-u_0(y))(\psi_\epsilon^-(x) -\psi_\epsilon^-(y))}{|x-y|^ {N+2s}}\,dxdy.
		    \end{aligned}
	        \end{equation}
			Let us define $\Omega_\epsilon := \{x\in\mathbb{R^N}:u_0^+(x) +\epsilon \phi(x)\leq 0\}$. On combining \eqref{eq3.13}, \eqref{eq3.15} and \eqref{eq3.12} we get
			\begin{align*}
				0 \leq & \|u_0\|^{2\theta}+\epsilon\|u_0\|^{2(\theta-1)}\langle u_0,\phi \rangle+\|u_0\|^{2(\theta-1)}\langle u_0,\psi_\epsilon^- \rangle-\lambda \int\limits_{\Omega} (u_0^+)^{-\gamma}(u_0^+(x) +\epsilon \phi(x)) \,dx\\
			&-\int\limits_{\Omega}\int\limits_{\Omega}\dfrac{|u_0^{+}(y)|^{2^{*}_{\mu ,s}}(u_0^+(x))^{2^{*}_{\mu ,s}-1}(u_0^+ +\epsilon \phi) (x)}{|x-y|^ \mu}\,dxdy + \lambda \int \limits_{\Omega_\epsilon} (u_0^+)^{-\gamma}(u_0^+(x) +\epsilon \phi(x)) \,dx\\ &+\int\limits_{\Omega}\int\limits_{\Omega_\epsilon} \dfrac{|u_0^{+}(y)|^{2^{*}_{\mu ,s}}(u_0^+(x))^{2^{*}_{\mu ,s}-1}(u_0^+ +\epsilon \phi) (x)}{|x-y|^ \mu}\,dxdy
			\end{align*}
				\begin{equation}\label{eq3.16}
					\begin{aligned}
						\quad\leq &\|u_0\|^{2\theta}-\lambda \int\limits_{\Omega} (u_0^+)^{1-\gamma}\,dx-\int\limits_{\Omega}\int\limits_{\Omega}\dfrac{|u_0^{+}(y)|^{2^{*}_{\mu ,s}}|u_0^+(x)|^{2^{*}_{\mu ,s}}}{|x-y|^ \mu}\,dxdy+\|u_0\|^{2(\theta-1)}\langle u_0,\psi_\epsilon^- \rangle\\
						&+\epsilon\left[\|u_0\|^{2(\theta-1)}\langle u_0,\phi \rangle-\lambda \int\limits_{\Omega}(u_0^+)^{-\gamma}\phi \,dx-\int\limits_{\Omega}\int\limits_{\Omega}\dfrac{|u_0^{+}(y)|^{2^{*}_{\mu ,s}}(u_0^+(x))^{2^{*}_{\mu ,s}-1}\phi(x)}{|x-y|^ \mu}\,dxdy\right] \\
					    \quad = &\|u_0\|^{2(\theta-1)}\langle u_0,\psi_\epsilon^- \rangle+\epsilon \left[  \|u_0\|^{2(\theta-1)}\langle u_0,\phi \rangle-\lambda \int\limits_{\Omega}(u_0^+)^{-\gamma}\phi \,dx\right] \\
						&-\epsilon\left[ \int\limits_{\Omega}\int\limits_{\Omega}\dfrac{|u_0^{+}(y)|^{2^{*}_{\mu ,s}}(u_0^+(x))^{2^{*}_{\mu ,s}-1}\phi(x)}{|x-y|^ \mu}\,dxdy \right].
					\end{aligned}
			        \end{equation}
				By the symmetry of the fractional kernel and following same steps as \eqref{eq3.14}
			    \begin{equation}\label{eq3.17}
			    	(u_0(x)- u_0(y))(u_0^+(x)-u_0^+(y))\geq {|u_0^+(x)-u_0^+(y)|}^2
			    \end{equation}
				Assuming,
				$$U(x,y)=\dfrac{(u_0(x)-u_0(y))(\phi(x)-\phi(y))}{|x-y|^{N+2s}}.$$
				Using \eqref{eq3.17} we get,
				\begin{equation}\label{eq3.18}
					\begin{aligned}
				\langle u_0,\psi_\epsilon^- \rangle =& \iint\limits_{\Omega_\epsilon\times\Omega_\epsilon}\dfrac{(u_0(x)-u_0(y))(\psi_\epsilon^-(x) -\psi_\epsilon^-(y))}{|x-y|^{N+2s}}\,dxdy\\&+2\iint\limits_{\Omega_\epsilon\times\mathbb{R^N}\backslash \Omega_\epsilon}\dfrac{(u_0(x)-u_0(y))(\psi_\epsilon^-(x) -\psi_\epsilon^-(y))}{|x-y|^{N+2s}}\,dxdy\\
				\leq&-\epsilon\iint\limits_{\Omega_\epsilon\times\Omega_\epsilon}U(x,y)\,dxdy-2\epsilon \iint\limits_{\Omega_\epsilon\times\mathbb{R^N}\backslash \Omega_\epsilon}U(x,y)\,dxdy\\
				\leq&\, 2\epsilon  \iint\limits_{\Omega_\epsilon\times\mathbb{R^N}}\left|U(x,y)\right| \,dxdy.
				\end{aligned}
				\end{equation}
		   	Now as $U \in L^1(\mathbb{R}^{2N})$, then for any $\sigma>0$ there exists $R_\sigma$ sufficiently large such that
			 $$\iint\limits_{(supp (\phi))\times(\mathbb{R^N}\backslash B_{R_{\sigma}})} |U(x,y)|\,dxdy< \frac{\sigma}{2}.$$
			  Since $\Omega_\epsilon \subset supp(\phi)$ and $|\Omega_\epsilon \times B_{R_{\sigma}}|\to 0$ as $\epsilon \to 0^+$.
			 Hence, there exist $\delta_\sigma > 0$ and $\epsilon_\sigma > 0$ such that for any $\epsilon \in (0,\epsilon_\sigma]$,
			$|\Omega_\epsilon \times B_{R_{\sigma}}|<\delta_\sigma$ and $\iint\limits_{\Omega_\epsilon\times B_{R_{\sigma}}} |U(x,y)|\,dxdy< \frac{\sigma}{2}.$
			Therefore, for any $\epsilon \in (0,\epsilon_\sigma]$,
			\begin{equation}\label{eq3.19}
				\begin{aligned}\iint\limits_{\Omega_\epsilon\times \mathbb{R^N}} |U(x,y)|\,dxdy < \sigma,\; \text{and}\;
					 \lim\limits_{\epsilon\to 0^+}\iint\limits_{\Omega_\epsilon\times\mathbb{R^N}} |U(x,y)|\,dxdy &=0.
					\end{aligned}
			\end{equation}
			 Therefore, dividing by $\epsilon$ and considering \eqref{eq3.16}, \eqref{eq3.18}, \eqref{eq3.19} with $\epsilon\to 0^+$ we get,
			\begin{equation*}
				0 \leq \|u_0\|^{2(\theta-1)}\langle u_0,\phi \rangle-\lambda \int\limits_{\Omega} (u_0^+)^{-\gamma}\phi(x) \,dx -\int\limits_{\Omega}\int\limits_{\Omega}\dfrac{|u_0^{+}(y)|^{2^{*}_{\mu ,s}}(u_0^+(x))^{2^{*}_{\mu ,s}-1}\phi(x)}{|x-y|^ \mu}\,dxdy,
			\end{equation*}
		for all $\phi \in X_0$.
		    As, $\phi$ was arbitrarily chosen, we infer that $u_0$ is the non trivial weak solution.
		    Finally taking $\phi= u_0^-$\;in \eqref{eq2.3} and using \eqref{eq3.14} 
		    \begin{align*}
		    	0 = \|u_0\|^{2(\theta-1)} \langle u_0, u_0^-\rangle &\leq -\|u_0\|^{2(\theta-1)}\|u_0^-\|^2.
		   \end{align*}
		    This implies $u_0 \geq 0$.
		     Moreover, by using the maximum principle in \cite[Proposition 2.17]{28prop2.17}, we get that $u_0$ is a positive solution of problem \eqref{eq2.2}.
		    \QED       
		    \section{Perturbed problem $(P_{n,\lambda})$}
		    In this section we prove that the functional $J_{n, \lambda}$ satisfies mountain-pass geometry and thus conclude about the existence of solution for the perturbed case.
		    \begin{Lemma}\label{Lemma4.1}
		    	Let $\rho \in(0, 1],\; \lambda_0 = \lambda_0(\rho) > 0$ and $\alpha=\alpha(\rho) > 0$ be given as in Lemma \ref{Lemma1}. Then, for any $\lambda\in (0, \lambda_0]$
		    	and any $u \in X_0$ with $\|u\| = \rho$, one has $J_{n,\lambda}(u) \geq \alpha$. Furthermore, there exists $e\in X_0$, with $\|e\| > \rho$, such that
		    	$J_{n,\lambda}(e) < 0.$
		    \end{Lemma}
		    
		    \proof For any $u\in X_0$ and any $n\in \mathbb{N}$, the subadditivity of $t \mapsto t^{1-\gamma}$ gives,
		    \begin{equation}\label{eq4.1}
		    	\displaystyle { {\left(u^+ +\frac{1}{n}\right)}^{1-\gamma}-\left(\frac{1}{n}\right)^{1-\gamma}\leq\left(u^+ \right)}^{1-\gamma}.
		    \end{equation}
		    Thus, we have
		    $$J_{n,\lambda}(u) \geq J_\lambda(u)\quad \text{for all}\; u \in X_0.$$
		    Therefore, by Lemma \ref{Lemma1} if $\|u\|=\rho$ we get	
		    $$J_{n,\lambda}(u)  \geq \alpha.$$
		    For every $v\in X_0$, with $v^+\neq 0$, and $t > 0$, we have
		    \begin{equation*}
		    	\displaystyle J_{n,\lambda}(tv)= {\frac{t^{2\theta}}{2\theta}\|v\|^ {2\theta}-\frac{\lambda}{1-\gamma}\int\limits_{\Omega}\left[ {\left({tv}^+ +\frac{1}{n}\right)}^{1-\gamma}-\left(\frac{1}{n}\right)^{1-\gamma}\right] \,dx -\frac{t^{2\cdot{2^*}_{\mu,s}}}{2\cdot{2^*}_{\mu,s}}\|v^{+}\|_{N L}^{2\cdot2^{*}_{\mu ,s}}} \to -\infty,
		    \end{equation*}
		    as $t \to \infty$. Therefore, we can find a sufficiently large $e\in X_0$, with $\|e\| > \rho$, such that $ J_{n,\lambda}(e)<0 $.\QED
		   Now we will establish two positive constants say $D_1$ and $D_2$ such that
		   \begin{equation}\label{eq4.3}
		   	\begin{aligned}
		   		{D_1} &=\displaystyle\left( \frac{1}{2\theta}-\frac{1}{2\cdot 2_{\mu,s}^*}\right)  {(S_s^H)}^\frac{2_{\mu,s}^*\theta}{2_{\mu,s}^*-\theta},\\
		   		{D_2} &=\displaystyle\frac{\left[ \left( \frac{1}{1-\gamma}+\frac{1}{2\cdot 2_{\mu,s}^*}\right)|\Omega|^\frac{2_s^*-1+\gamma}{2_s^*}S_s^\frac{-(1-\gamma)}{2}\right]^\frac{2\theta}{2\theta-1+\gamma}}{\left( \frac{1}{2\theta}-\frac{1}{2\cdot 2_{\mu,s}^*}\right)^\frac{1-\gamma}{2\theta-1+\gamma}}.
		   	\end{aligned}
		   \end{equation}
		   
		   \begin{Lemma}\label{Lemma4.2}
		   	Let $\lambda > 0$. Then the functional $J_{n,\lambda}$ satisfies the Palais–Smale condition at any level $c \in \mathbb{R}$ verifying
		   	\begin{equation}\label{eq4.4}
		   		c < D_1-D_2 \lambda^\frac{2\theta}{2\theta-1+\gamma},
		   	\end{equation} with $D_1, D_2 > 0$ given as in \eqref{eq4.3}.
		   \end{Lemma}
		   
		   \proof:   Let ${(u_k)}_k$ be a Palais–Smale sequence in $X_0$ at level $c\in \mathbb{R}$, satisfying \eqref{eq4.4}. We have to prove that ${(u_k)}_k$ has a convergent sub-sequence in $X_0$.
		   So we first claim the boundedness of ${(u_k)}_k$, by using \eqref{eq2.7}, \eqref{eq4.1}, the H\"older inequality and by \eqref{eq4.2}, we say there exists $\sigma > 0$ such that, as $k\to \infty$,
		   \begin{align*}
		   	c+\sigma\|u_k\|+o(1)\geq& J_{n,\lambda }(u_k)-\frac{1}{2\cdot{2^*}_{\mu,s}}\langle J_{n,\lambda }^{'}(u_k),(u_k) \rangle\\
		   	%
		   	%
		   	\geq& \left(\frac{1}{2\theta}-\frac{1}{2\cdot 2_{\mu,s}^*}\right)\|u_k\|^{2\theta}-\lambda\left( \frac{1}{1-\gamma}+\frac{1}{2\cdot 2_{\mu,s}^*}\right)|\Omega|^\frac{2_s^*-1+\gamma}{2_s^*}S_s^\frac{-(1-\gamma)}{2}\|u_k\|^{1-\gamma}.
		   \end{align*}
		  Therefore, ${(u_k)}_k$ is bounded in $X_0$, since $1-\gamma<1<2\theta$.
		   Besides $(u_k)^-$ is also bounded in $X_0$, so by \eqref{eq4.2} we have,
		   \begin{equation*}
		   	\lim\limits_{k \to\infty}\langle J_{n,\lambda }^{'}(u_k),-u_k^- \rangle=\lim\limits_{k \to\infty}\|u_k\|^{2\theta-1}\langle u_k,-u_k^-\rangle=0.
		   \end{equation*}
		    As $k \to \infty$ we get $\|u_k^-\|\to 0$ by inequality \eqref{eq3.14}, so we conclude
		   \begin{align*}
		   	J_{n,\lambda }(u_k) =J_{n,\lambda }(u_k^+)+o(1),\;\;
		   	J'_{n,\lambda }(u_k) =J'_{n,\lambda }(u_k^+)+o(1)\quad \text{as} \; k\to \infty.
		   \end{align*}
		   Therefore, we can assume $(u_k)$ to be a sequence of non negative functions which is bounded. Hence, up to a sub-sequence, there exists a function $u\in X_0$ such that
		   \begin{equation}\label{eq4.5}\left\{
		   	\begin{array}{lr}
		   		u_k\rightharpoonup u\; \text{in} \; X_0,\; \|u_k\|\rightarrow \mu,\;
		   		u_k\rightharpoonup u\; \text{in} \; L^{2_s^*}(\Omega),\; u_k\rightarrow u\; \text{in} \; L^{p}(\Omega)\; \text{for all}\; p\to [1,2_s^*),\\
		   		\displaystyle{\int\limits_{\Omega} \frac{|u_k(y)|^{2^{*}_{\mu ,s}}|u_k(x)|^{2^{*}_{\mu ,s}-2}u_k(x)}{|x-y|^ \mu}\,dy \rightharpoonup \int\limits_{\Omega} \frac{|u(y)|^{2^{*}_{\mu ,s}}|u(x)|^{2^{*}_{\mu ,s}-2}u(x)}{|x-y|^ \mu}\,dy}\; \text{weakly in} \; L^\frac{2N}{N+2s},\\
		   		\|u_k-u\|_{N L}\rightarrow d ,\\
		   		u_k\rightarrow u \quad\text{a.e} \;\; \text{and}\;\; u_k \leq h \;\;\text{a.e in} \; \Omega \;\;\text{as} \; k\to \infty,\; \text{with}\; h\in L^p(\Omega).
		   	\end{array}
		   	\right.
		   \end{equation}
		   If $\mu =0$, then $u_k\rightarrow 0$ in $X_0$.
		   Let us assume that $\mu> 0$. Using \eqref{eq4.5} we get
		   \begin{align*}
		   	\displaystyle{ \left|\frac{u_k-u}{\left(u_k +\frac{1}{n}\right)^{\gamma}}\right|\leq n^\gamma|u_k-u|\leq n^\gamma(h+|u|) }\quad\text{a.e in} \; \Omega,
		   \end{align*}
		   by making use of the dominated convergence theorem and \eqref{eq4.5} we get,
		   \begin{equation}\label{eq4.6}
		   	\lim\limits_{k \to\infty}\int\limits_{\Omega}\frac{u_k-u}{\left(u_k +\frac{1}{n}\right)^{\gamma}}\,dx=0.
		   \end{equation}
		   By Br\'{e}zis Lieb lemma \cite[Theorem2]{7thm2}, \cite[Lemma2.2]{M.yang} and \eqref{eq4.5} as $k \to \infty$,
		   \begin{equation}\label{eq4.7}
		   		\displaystyle\|u_k\|^2 = \|u_k-u\|^2+\|u\|^2+ o(1),\;
		   		\displaystyle \|u_k\|_{N L}^{2\cdot2^{*}_{\mu ,s}} = \|u_k-u\|_{N L}^{2\cdot2^{*}_{\mu ,s}} + \|u\|_{N L}^{2\cdot2^{*}_{\mu ,s}} +o(1).
		   \end{equation}
		   From \eqref{eq4.2}, \eqref{eq4.5}, \eqref{eq4.6} and \eqref{eq4.7} we derive that, as $k \to \infty$
		   \begin{align*}
		   	o(1)  = &\displaystyle \langle J_{n,\lambda }^{'}(u_k),u_k-u \rangle\\
		   	 = & \displaystyle\|u_k\|^{2(\theta-1)}\langle u_k,u_k-u\rangle-\lambda\int\limits_{\Omega}\frac{u_k-u}{\left(u_k +\frac{1}{n}\right)^{\gamma}}\,dx\\
		   	&-\int\limits_{\Omega}\int\limits_{\Omega}\dfrac{|u_k(y)|^{2^{*}_{\mu ,s}}(u_k(x))^{2^{*}_{\mu ,s}-1}(u_k(x)-u(x))}{|x-y|^ \mu}\,dxdy\\
		   	 = &\displaystyle\mu^{2(\theta-1)}[\mu^{2}-\|u\|^{2}]-\|u_k\|_{N L}^{2\cdot2^{*}_{\mu ,s}}+\int\limits_{\Omega}\int\limits_{\Omega}\dfrac{|u_k(y)|^{2^{*}_{\mu ,s}}(u_k(x))^{2^{*}_{\mu ,s}-1}u(x)}{|x-y|^ \mu}\,dxdy\\
		   	 = & \displaystyle\mu^{2(\theta-1)}[\mu^{2}-\|u\|^{2}] - \|u_k\|_{N,L}^{2\cdot2^{*}_{\mu ,s}} + \|u\|_{N L}^{2\cdot2^{*}_{\mu ,s}}+o(1)\\
		   	 = & \displaystyle\mu^{2(\theta-1)}\|u_k-u\|^2 - \|u_k-u\|_{N L}^{2\cdot2^{*}_{\mu ,s}}+o(1).
		   \end{align*}
		   From there we get a relation,
		   \begin{equation}\label{eq4.8}
		   	\displaystyle{\mu^{2(\theta-1)}\lim\limits_{k \to\infty}\|u_k-u\|^2 = \lim\limits_{k \to\infty}\|u_k-u\|_{N L}^{2\cdot 2^{*}_{\mu ,s}}} = d^{2\cdot 2^{*}_{\mu ,s}}.
		   \end{equation}
		   If $d=0$, then by \eqref{eq4.5}, \eqref{eq4.8} and positivity of $\mu$ it follows that
		   $\lim\limits_{k \to\infty}\|u_k-u\|^2=0$,\;i.e. $u_k\rightarrow u$ in $X_0$ which concludes the proof.
		   So let us assume by contradiction that $d > 0$, then using \eqref{eq2.7} we get
		   \begin{equation}\label{eq4.9}
		   	\displaystyle d^{2\cdot 2^{*}_{\mu ,s}} \geq S_s^H \mu^{2(\theta-1)}d^{2},
		   \end{equation}
		   i.e.
		   \begin{equation*}
		   	\displaystyle d^{2\cdot (2^{*}_{\mu ,s}-1)} \geq S_s^H \mu^{2(\theta-1)},
		   \end{equation*}
		   also \eqref{eq4.8} implies that
		   \begin{equation*}
		   	\displaystyle d^{2\cdot 2^{*}_{\mu ,s}}=\mu^{2(\theta-1)}[\mu^{2}-\|u\|^{2}],
		   \end{equation*}
		   from \eqref{eq4.9}, we have
		   \begin{equation}\label{eq4.10}
		   	\displaystyle \mu^{2}\geq (S_s^H)^\frac{2^{*}_{\mu ,s}}{2^{*}_{\mu ,s}-\theta}.
		   \end{equation}
		   As $1<\theta < 2^{*}_{\mu ,s}$ , we get $2^{*}_{\mu ,s}-\theta>0$
		   \begin{align*}
		   	\displaystyle J_{n,\lambda }(u_k)-\frac{1}{2\cdot{2^*}_{\mu,s}}\langle J_{n,\lambda }^{'}(u_k),(u_k) \rangle\\
		   		\geq \left(\frac{1}{2\theta}-\frac{1}{2\cdot 2_{\mu,s}^*}\right)\|u_k\|^{2\theta}&-\lambda\left( \frac{1}{1-\gamma}+\frac{1}{2\cdot 2_{\mu,s}^*}\right)\int\limits_{\Omega}|u_k|^{1-\gamma}\,dx.
		   \end{align*}
		   By \eqref{eq4.2}, \eqref{eq4.5}, \eqref{eq4.7}, \eqref{eq4.10} and the Young inequality with $k \to \infty$ we get,
		   \begin{align*}
		   	c \geq & \displaystyle\left(\frac{1}{2\theta}-\frac{1}{2\cdot 2_{\mu,s}^*}\right)(\mu^{2\theta}+\|u\|^{2\theta})-\lambda\left( \frac{1}{1-\gamma}+\frac{1}{2\cdot 2_{\mu,s}^*}\right)|\Omega|^\frac{2_s^*-1+\gamma}{2_s^*}S_s^\frac{-(1-\gamma)}{2}\|u\|^{1-\gamma}\\
		   	 \geq & \displaystyle\left(\frac{1}{2\theta}-\frac{1}{2\cdot 2_{\mu,s}^*}\right)(\mu^{2\theta}+\|u\|^{2\theta})-\left(\frac{1}{2\theta}-\frac{1}{2\cdot 2_{\mu,s}^*}\right)\|u\|^{2\theta}\\
		   	 & -\displaystyle{\left(\frac{1}{2\theta}-\frac{1}{2\cdot 2_{\mu,s}^*}\right)}^\frac{-(1-\gamma)}{2\theta-1+\gamma}{{\left(\lambda\left( \frac{1}{1-\gamma}+\frac{1}{2\cdot 2_{\mu,s}^*}\right)|\Omega|^\frac{2_s^*-1+\gamma}{2_s^*}S_s^\frac{-(1-\gamma)}{2}\right)}^\frac{2\theta}{2\theta-1+\gamma}}\\
		      \displaystyle\geq & \left(\frac{1}{2\theta}-\frac{1}{2\cdot 2_{\mu,s}^*}\right){(S_s^H)}^\frac{2_{\mu,s}^*\theta}{2_{\mu,s}^*-\theta}\\&-{\left(\frac{1}{2\theta}-\frac{1}{2\cdot 2_{\mu,s}^*}\right)}^\frac{-(1-\gamma)}{2\theta-1+\gamma}{{\left[\lambda\left( \frac{1}{1-\gamma}+\frac{1}{2\cdot 2_{\mu,s}^*}\right)|\Omega|^\frac{2_s^*-1+\gamma}{2_s^*}S_s^\frac{-(1-\gamma)}{2}\right]}^\frac{2\theta}{2\theta-1+\gamma}}\\
		   	 = &\;D_1- \lambda^\frac{2\theta}{2\theta-1+\gamma} D_2
		   	 >  c.
		   \end{align*}
		   which is a contradiction and thus concludes the proof.\QED

		  \begin{Lemma}\label{Lemma4.3}There exist $\Psi\in X_0, \lambda^{*}>0,$ such that $\lambda\in(0,\lambda^{*})$
		  	$$\sup\limits_{t\geq 0} J_{n,\lambda}(t\Psi)< D_1-D_2\lambda^\frac{2\theta}{2\theta-1+\gamma},$$
		  	with $D_1, D_2$ given as in \eqref{eq4.3}.
		  	\end{Lemma}
		  
		  \proof: $u_\epsilon$ be as defined in \eqref{eq2.9} and from \eqref{eq2.6}, 
		  \begin{align*}
		  	\displaystyle J_{n,\lambda}(tu_\epsilon) =&\; \frac{t^{2\theta}}{2\theta}\|u_\epsilon\|^ {2\theta}-\frac{\lambda}{1-\gamma}\int\limits_{\Omega}\left[{\left(tu_\epsilon +\frac{1}{n}\right)}^{1-\gamma}-\left(\frac{1}{n}\right)^{1-\gamma}\right]dx\\ &-\frac{t^{2\cdot{2_{\mu,s}^*}}}{2\cdot{2_{\mu,s}^*}}\int\limits_{\Omega}\int\limits_{\Omega}\frac{|u_\epsilon(y)|^{2_{\mu,s}^*}|u_\epsilon(x)|^{2_{\mu,s}^*}}{|x-y|^\mu}\,dxdy\\
		  	< &\; ct^{2\theta}.
		  \end{align*}
		  We can choose $\lambda^1 > 0$ such that $ 0< D_1-D_2\lambda^\frac{2\theta}{2\theta-1+\gamma}$\, and $t_0\in (0,1)$ such that\\$\displaystyle\sup\limits_{0\leq t\leq t_0} J_{n,\lambda}(tu_\epsilon) < D_1-D_2\lambda^\frac{2\theta}{2\theta-1+\gamma}$ for all $\lambda \in (0, \lambda^1)$. This implies we only have to show $\displaystyle {\sup\limits_{t\geq t_0} J_{n,\lambda}(tu_\epsilon)< D_1-D_2\lambda^\frac{2\theta}{2\theta-1+\gamma}}$.
		  \begin{align*}
		  \sup\limits_{t\geq t_0}J_{n,\lambda}(tu_\epsilon)  &= \sup\limits_{t\geq t_0}\left( \frac{t^{2\theta}}{2\theta}\|u_\epsilon\|^ {2\theta}-\frac{\lambda}{1-\gamma}\int\limits_{\Omega} {\left(tu_\epsilon +\frac{1}{n}\right)}^{1-\gamma}-\left(\frac{1}{n}\right)^{1-\gamma}dx -\frac{t^{2\cdot{2_{\mu,s}^*}}}{2\cdot{2_{\mu,s}^*}}\|u_\epsilon\|_{N L}^{2\cdot{2_{\mu,s}^*}}\right) \\
		 & \leq \displaystyle\sup\limits_{t\geq 0}\nu(t)-\frac{\lambda}{1-\gamma}\int\limits_{\Omega}\left[ {\left(t_0u_\epsilon +\frac{1}{n}\right)}^{1-\gamma}-\left(\frac{1}{n}\right)^{1-\gamma}\right] \,dx,
		  \end{align*}
		  where
		  \begin{equation*}
		  	\nu(t)=\displaystyle\frac{t^{2\theta}}{2\theta}\|u_\epsilon\|^ {2\theta}-\frac{t^{2\cdot{2_{\mu,s}^*}}}{2\cdot{2_{\mu,s}^*}}\|u_\epsilon\|_{N L}^{2\cdot{2_{\mu,s}^*}}.
		  \end{equation*}
		  Since $\nu(0)=0$, $\nu(t)\to -\infty$ as $t\to\infty$ and $\nu(t)>0$ close to zero.
		  Therefore, there exists $t_\epsilon>0$ such that $\sup\limits_{t\geq 0}\nu(t)=\nu(t_\epsilon)$. Accordingly $\nu'(t_\epsilon)=0$, which gives
		 \begin{equation}\label{eq4.11}
		 \displaystyle t_\epsilon  = {\left[\frac{\|u_\epsilon\|^ {2\theta}}{\|u_\epsilon\|_{N L}^{2\cdot{2_{\mu,s}^*}}}\right]}^\frac{1}{2\cdot{2_{\mu,s}^*}-2\theta}, \quad \;
		 \displaystyle\nu(t_\epsilon)  = \left(\frac{1}{2\theta}-\frac{1}{2\cdot 2_{\mu,s}^*}\right){\left[\frac{\|u_\epsilon\|^ {2}}{\|u_\epsilon\|_{N L}^{2}}\right]}^\frac{{2_{\mu,s}^*}\theta}{{2_{\mu,s}^*}-\theta},
		 \end{equation}
		  Also from the Proposition \ref{Prop2.3},
		  \[\|u_\epsilon\|^ {2}\leq C^\frac{N(N-2s)}{2s(2N-\mu)}(S^H_s)^\frac{N}{2s}+ O(\epsilon^{N-2s}),
		  \;\; \|u_\epsilon\|_{N L}^{2}\geq {\left[C^\frac{N}{2s}(S^H_s)^\frac{2N-\mu}{2s}- O(\epsilon^{N})\right] }^\frac{N-2s}{2N-\mu}.\]
		  This implies,
		 \begin{align*}
		 	\frac{\|u_\epsilon\|^ {2}}{\|u_\epsilon\|_{N L}^{2}}\leq \frac{C^\frac{N(N-2s)}{2s(2N-\mu)}(S^H_s)^\frac{N}{2s}+ O(\epsilon^{N-2s})}{{\left[C^\frac{N}{2s}(S^H_s)^\frac{2N-\mu}{2s}- O(\epsilon^{N})\right] }^\frac{N-2s}{2N-\mu}}
		 	= S^H_s \left[1+ O(\epsilon^{N-2s})\right],
		 \end{align*}		 
		  as  $\epsilon\in(0,1)$  and $ N-2s< N.$ Therefore, there exists a positive constant $c_1$, such that
		  from \eqref{eq4.11} we get
		  \begin{equation*}
		  	\nu(t_\epsilon)\leq\left(\frac{1}{2\theta}-\frac{1}{2\cdot 2_{\mu,s}^*}\right)(S^H_s)^\frac{{2_{\mu,s}^*}\cdot\theta}{{2_{\mu,s}^*}-\theta} + c_1\epsilon^{N-2s}.
		  \end{equation*}
		  Thus, we have for $t \geq t_0$
		  \begin{equation}\label{eq4.12}
		  	\begin{aligned}
		  		\displaystyle J_{n,\lambda}(tu_\epsilon) & \leq\, \sup\limits_{t\geq 0}\nu(t)-\frac{\lambda}{1-\gamma}\int\limits_{\Omega}\left[ {\left(t_0 u_\epsilon +\frac{1}{n}\right)}^{1-\gamma}-\left(\frac{1}{n}\right)^{1-\gamma}\right] \,dx \\
		  		\displaystyle & \leq\left(\frac{1}{2\theta}-\frac{1}{2\cdot 2_{\mu,s}^*}\right)(S^H_s)^\frac{{2_{\mu,s}^*}\cdot\theta}{{2_{\mu,s}^*}-\theta} + c_1\epsilon^{N-2s}{-\frac{\lambda}{1-\gamma}\int\limits_{\Omega}{\left(t_0 u_\epsilon +\frac{1}{n}\right)}^{1-\gamma}-\left(\frac{1}{n}\right)^{1-\gamma}dx}.
		  	\end{aligned}
		  \end{equation}
		  For any $p>1$, $a$, $b>0$ large enough, using inequality
		  \begin{equation*}
		  	a^{1-\gamma}-(a+b)^{1-\gamma}\leq -(1-\gamma)b^\frac{1-\gamma}{p}a^\frac{(p-1)(1-\gamma)}{p},
		  \end{equation*}
		  and taking $a=\frac{1}{n}$ , $b= t_0u_\epsilon$, and $p={2_{\mu,s}^*}$, we have
		  \begin{equation}\label{eq4.13}
		  	\begin{aligned}
		  		\displaystyle-\frac{\lambda}{1-\gamma}\int\limits_{\Omega}\left[ {\left(t_0 u_\epsilon +\frac{1}{n}\right)}^{1-\gamma}-\left(\frac{1}{n}\right)^{1-\gamma}\right] \,dx &\leq -\lambda C\int\limits_{\Omega}|u_\epsilon|^\frac{1-\gamma}{{2_{\mu,s}^*}}\,dx\\
		  		& \displaystyle\leq-\lambda C\int\limits_{B_{\delta}(0)}| U_\epsilon|^\frac{1-\gamma}{{2_{\mu,s}^*}}\,dx\\
		  		& = \displaystyle-\lambda C\epsilon^\frac{(N-2s)(1-\gamma)}{2\cdot{2_{\mu,s}^*}}\int\limits_{|x|\leq\delta}\frac{dx}{[\epsilon^{2}+|x|^{2}]^\frac{(N-2s)(1-\gamma)}{2\cdot{2_{\mu,s}^*}}}.
		  	\end{aligned}
		  \end{equation}
		  Considering, $\epsilon< \delta^\frac{1}{q}$ sufficiently small, where the positive number $q<1$ and since $0<1-\gamma<1<\theta<{2_{\mu,s}^*}$, $N>2s$ so $q$ also satisfies,
		  $$1+\frac{2\theta\left[ (N-2s)(1-\gamma)-2q(N-2s)(1-\gamma)+2\cdot{2_{\mu,s}^*}qN\right]}{(2\theta-1+\gamma)(N-2s)2\cdot{2_{\mu,s}^*}} -\frac{2\theta}{2\theta-1+\gamma}<0.$$
		  That is, $$\displaystyle 0< q< \min\left\lbrace 1,\,  \frac{(N-2s)({2^*}_{\mu,s}-\theta)(1-\gamma)}{2\theta[N({2_{\mu,s}^*}-1)+N\gamma+2s(1-\gamma)]}\right\rbrace. $$
		  From \eqref{eq4.13} we have,
		  \begin{align*}
		  	\displaystyle-\frac{\lambda}{1-\gamma}\int\limits_{\Omega}{\left(t_0 u_\epsilon +\frac{1}{n}\right)}^{1-\gamma}-\left(\frac{1}{n}\right)^{1-\gamma}\,dx &\leq-\lambda C\epsilon^\frac{(N-2s)(1-\gamma)}{2\cdot{2_{\mu,s}^*}}\int\limits_{\{x\in \Omega;|x|\leq\epsilon^q\}}\frac{dx}{[\epsilon^{2}+|x|^{2}]^\frac{(N-2s)(1-\gamma)}{2\cdot{2_{\mu,s}^*}}} \\
		  	& \leq\displaystyle -\lambda c_2 \epsilon^\frac{(N-2s)(1-\gamma)-2q(N-2s)(1-\gamma)+2\cdot{2_{\mu,s}^*}qN}{2\cdot{2_{\mu,s}^*}},
		  \end{align*}
		  where $c_2$ is a positive constant independent of $\epsilon$. Substituting these values in \eqref{eq4.12}, for $t \geq t_0$
		  \begin{equation*}
		  	\displaystyle J_{n,\lambda}(tu_\epsilon) \leq\left(\frac{1}{2\theta}-\frac{1}{2\cdot 2_{\mu,s}^*}\right)(S_s^H)^\frac{{2_{\mu,s}^*}\cdot\theta}{{2_{\mu,s}^*}-\theta} + c_1\epsilon^{(N-2s)}-\lambda c_2 \epsilon^\frac{(N-2s)(1-\gamma)-2q(N-2s)(1-\gamma)+2\cdot{2_{\mu,s}^*}qN}{2\cdot{2_{\mu,s}^*}},
		  \end{equation*}
		  and taking,\begin{equation*}
		  	\displaystyle {\epsilon= \lambda^\frac{2\theta}{(2\theta-1+\gamma)(N-2s)}},\quad
		  	\nu_2=1+\frac{2\theta\left[ (N-2s)(1-\gamma)-2q(N-2s)(1-\gamma)+2\cdot{2_{\mu,s}^*}qN\right]}{(2\theta-1+\gamma)(N-2s)2\cdot{2_{\mu,s}^*}} ,
		  \end{equation*}
		  \begin{equation*}
		  	\nu_3= \nu_2-\frac{2\theta}{(2\theta-1+\gamma)}.
		  \end{equation*}
		  And by assumption of q, we have $\nu_3 < 0$. Combining all these for the case $t \geq t_0,$ we get
		  \begin{align*}
		  	\displaystyle J_{n,\lambda}(tu_\epsilon) &\leq\left(\frac{1}{2\theta}-\frac{1}{2\cdot 2_{\mu,s}^*}\right)(S^H_s)^\frac{{2_{\mu,s}^*}\cdot\theta}{{2_{\mu,s}^*}-\theta} +c_1\lambda^\frac{2\theta}{(2\theta-1+\gamma)}-c_2\lambda^{\nu_2}\\
		  	\displaystyle J_{n,\lambda}(tu_\epsilon) & \leq D_1 +c_1\lambda^\frac{2\theta}{(2\theta-1+\gamma)}-c_2\lambda^{\nu_2}\\
		  	\displaystyle & = D_1 + \lambda^\frac{2\theta}{(2\theta-1+\gamma)}[-c_2\lambda^{\nu_2-\frac{2\theta}{(2\theta-1+\gamma)}}+c_1]\\
		  	\displaystyle & = D_1 - \lambda^\frac{2\theta}{(2\theta-1+\gamma)}[c_2\lambda^{\nu_3}-c_1].
		  \end{align*}
		  Thus, we choose $\lambda^{2}>0 $ such that for all $ \lambda\in(0,\lambda^2),$
		  \begin{equation*}
		  \begin{aligned}
		  		J_{n,\lambda}(tu_\epsilon) &< D_1 - \lambda^\frac{2\theta}{(2\theta-1+\gamma)}D_2,\quad \text{for all}\; t \geq t_0,\\
		  	\sup\limits_{t\geq t_0} J_{n,\lambda}(tu)&< D_1-\lambda^\frac{2\theta}{2\theta-1+\gamma}D_2.
		  \end{aligned}
		  \end{equation*}
	  Take $\lambda^{*} = \min\{\lambda^1, \lambda^2\}$,
	  \begin{equation*}
	  	\sup\limits_{t\geq 0} J_{n,\lambda}(tu) < D_1-\lambda^\frac{2\theta}{2\theta-1+\gamma}D_2, \quad \text{for all}\; \lambda \in (0,\lambda^{*}).
	  \end{equation*}\QED	
		\begin{Theorem}\label{Thm4.4}
			There exists $\overline{\lambda}>0$ such that for all $\lambda\in(0,\overline{\lambda}),$ problem \eqref{pb2.4} has a positive solution $v_n\in X_0$ with 
			\begin{equation}\label{eq4.14}
				\alpha< J_{n,\lambda}{(v_n)}< D_1-\lambda^\frac{2\theta}{2\theta-1+\gamma}D_2,
			\end{equation}
			where $\alpha,$ $D_1$ and $D_2$ are given in Lemma \ref{Lemma1} and \eqref{eq4.3} respectively.
		\end{Theorem} 
		
		\proof:  Clearly from Lemma \ref{Lemma4.1}, we can see $J_{n,\lambda}$ satisfies the geometry of the mountain-pass lemma. Set the critical level as,
		$$c_{n,\lambda}=\inf\limits_{g\in\Gamma} \max \limits_{t \in [0,1]} J_{n,\lambda}(g(t)). $$
		Here $\Gamma= \{g\in C([0,1],X_0): g(0)=0,\, J_{n,\lambda}(g(1))<0 \}$ and $\lambda\in(0,\overline{\lambda})$ with $\overline{\lambda}= \min \{\lambda_0,\lambda^{*} \}$ where $\lambda_0$ and $\lambda^*$ are as given in Lemma \ref{Lemma1} and \ref{Lemma4.3}. Besides by using Lemma \ref{Lemma4.1} and \ref{Lemma4.3} we claim,
		$$0< \alpha<c_{n,\lambda}\leq \sup\limits_{t\geq 0} J_{n,\lambda}(t\Psi)< D_1-D_2\lambda^\frac{2\theta}{2\theta-1+\gamma}. $$
		Hence, by Lemma \ref{Lemma4.2} we conclude that $J_{n,\lambda}$ satisfies Palais -Smale condition at level $c_{n,\lambda}.$ Thus, by applying the mountain-pass lemma, we get a critical point $v_n \in X_0 $ for $J_{n,\lambda}$ at level $c_{n,\lambda}$
		i.e. \begin{equation*}
			J'_{n,\lambda}(v_n)=0 ,\quad J_{n,\lambda}(v_n)=c_{n,\lambda}>\alpha>0=J_{n,\lambda}(0).
		\end{equation*}
		Therefore, $v_n$ is non-trivial solution of problem \eqref{pb2.4}.
		Taking the test function $\phi = v_n^-$ and by
		inequality \eqref{eq3.14}, we get $\|v_n^-\|=0$, hence $v_n \geq 0$. Finally using the maximum principle \cite[Proposition2.17]{28prop2.17}, we obtain that $v_n$ is a positive solution of \eqref{pb2.4}.\QED
		\section{Existence of  Second Solution}
		In this section we show the existence of second solution of problem $(P_\lambda)$, and this we do by passing the limit in the family of solutions of problem \eqref{pb2.4}.
		
		{\bf Proof of  Theorem \ref{thm1.1}:} Let $v_n$ be the family of positive solutions of problem \eqref{pb2.4} and let $\lambda \in (0,\overline\lambda)$ where $\overline\lambda$ be as given in Theorem \ref{Thm4.4}. By \eqref{eq2.8}, \eqref{eq4.1}, \eqref{eq4.14} and the H\"older inequality we have,
		\begin{align*}
			D_1-D_2\lambda^\frac{2\theta}{2\theta-1+\gamma} >&J_{n,\lambda}{(v_n)}-\frac{1}{2\cdot{2^*}_{\mu,s}}\langle J_{n,\lambda }^{'}(v_n),(v_n) \rangle\\
			 \geq &\left(\frac{1}{2\theta}-\frac{1}{2\cdot 2_{\mu,s}^*}\right)\|v_n\|^{2\theta}-\lambda\left( \frac{1}{1-\gamma}+\frac{1}{2\cdot 2_{\mu,s}^*}\right)|\Omega|^\frac{2_s^*-1+\gamma}{2_s^*}S_s^\frac{-(1-\gamma)}{2}\|v_n\|^{1-\gamma}.
		\end{align*}
		We see that ${(v_n)}_n$ is bounded in $X_0$ and up to a sub-sequence, there exists a function $v_0\in X_0$ such that 
		\begin{equation}\label{eq5.1}\left\{
			\begin{array}{lr}
				v_n\rightharpoonup v_0\; \text{in} \; X_0,\;\; \|v_n\|\rightarrow \beta,\;\; v_n\rightharpoonup v_0 \; \text{in} \; L^{2_s^*}(\Omega),\;\; v_n\rightarrow v_0\; \text{in} \; L^{p}(\Omega)\; \text{for all}\; p\in [1,2_s^*),\\
				
				\displaystyle{\int\limits_{\Omega} \dfrac{|v_n(y)|^{2^{*}_{\mu ,s}}|v_n(x)|^{2^{*}_{\mu ,s}-2}v_n(x)}{|x-y|^ \mu}\,dy \rightharpoonup \int\limits_{\Omega} \dfrac{|v_0(y)|^{2^{*}_{\mu ,s}}|v_0(x)|^{2^{*}_{\mu ,s}-2}v_0(x)}{|x-y|^ \mu}\,dy}\; \text{weakly in}\, L^\frac{2N}{N+2s},\\
				\|v_n-v_0\|_{N L}\rightarrow d,\\
				
			\end{array}
			\right.
		\end{equation}
		If $\beta =0$, then $v_n\rightarrow 0$ in $X_0$ as $n\to\infty$. So, let us assume that $\beta> 0$.
		We note that $\displaystyle 0\leq\frac{v_n}{(v_n +\frac{1}{n})^\gamma}\leq v_n^{1-\gamma} \; a.e \; in \; \Omega,$
		and by the Vitali convergence theorem we have,
		\begin{equation}\label{eq5.2}
			\lim\limits_{n\to \infty}\int\limits_{\Omega}\frac{v_n}{(v_n +\frac{1}{n})^\gamma}\,dx=\int\limits_{\Omega} v_0^{1-\gamma}.
		\end{equation}
		Considering \eqref{eq2.5} for $v_n$ and taking test function $\phi = v_n$ we get,
		$$ {\|v_n\|^{2\theta}-\lambda \int\limits_{\Omega} \left(v_n +\frac{1}{n}\right)^{-\gamma}v_n \,dx -\int\limits_{\Omega}\int\limits_{\Omega}  \dfrac{|v_n(y)|^{2^{*}_{\mu ,s}}|v_n(x)|^{2^{*}_{\mu ,s}}}{|x-y|^ \mu}\,dxdy=0},$$
		by \eqref{eq5.1}, \eqref{eq5.2} as $n\to \infty$,
		\begin{equation}\label{eq5.3}
			\beta^{2\theta}-\lambda \int\limits_{\Omega} (v_0)^{1-\gamma}\,dx -\|v_n\|_{N L}^{2\cdot2^*_{\mu,s}}= o(1).
		\end{equation}
		With an easy computation we see that for, any $\displaystyle n\in \mathbb{N}$, there exists a positive constant C such that
		\begin{equation*}
			\displaystyle \|v_n\|^{2(\theta-1)}(-\Delta)^sv_n\geq \min\left\lbrace C, \frac{\lambda}{2^\gamma}\right\rbrace\; in\; \Omega.
		\end{equation*}
		Therefore, by standard comparison argument and maximum principle result \cite{3Lemma2.1, 28prop2.17} for any $\tilde{\Omega}\Subset \Omega$ there exists a constant $c_{\tilde{\Omega}}> 0$ such that
		\begin{equation}\label{eq5.4}
			v_n\geq c_{\tilde{\Omega}}\quad a.e.\; \text{in} \; \tilde{\Omega}\; \text{and for any}\; n \in \mathbb{N}.
		\end{equation}
		Using \eqref{eq5.4}, for $\phi \in C_0^{\infty}(\Omega)$ with $supp\, \phi = \tilde{\Omega}\Subset \Omega$ we get $\displaystyle 0\leq\left| \dfrac{\phi}{(v_n +\frac{1}{n})^\gamma}\right| \leq \dfrac{|\phi|}{{c^{\gamma}_{\tilde{\Omega}}}},$ so by the Dominated convergence theorem and \eqref{eq5.1} we get,
		\begin{equation}\label{eq5.5}
			\lim\limits_{n\to \infty}\int\limits_{\Omega}\frac{\phi}{(v_n +\frac{1}{n})^\gamma}\,dx = \int\limits_{\Omega}
			\frac{\phi}{{v_0}^{\gamma}}\,dx.
		\end{equation}
	Passing the limit $n \to \infty$ in \eqref{eq2.5} for $v_n$ and from  \eqref{eq5.1} and \eqref{eq5.5} one gets,
		\begin{equation}\label{eq5.6}
			\displaystyle {\beta^{2(\theta-1)}\langle v_0, \phi \rangle-\lambda \int\limits_{\Omega} v_ 0^{-\gamma}\phi(x) \,dx -\int\limits_{\Omega}\int\limits_{\Omega}  \dfrac{|v_0(y)|^{2^{*}_{\mu ,s}}(v_0(x))^{2^{*}_{\mu ,s}-1}\phi(x)}{|x-y|^ \mu}\,dxdy} = 0.
		\end{equation}
	By the standard density argument it holds for all $\phi \in X_0$, since boundary is in $ C^2$ and that implies $C_0^{\infty}(\Omega)$ is dense in $X_0$ \cite{15Thm6}. Thus in particular it also holds for $\phi=v_0$. Passing limit $n \to \infty$ and combining \eqref{eq5.3} and \eqref{eq5.6} we get,
		\begin{equation*}
			\beta^{2(\theta-1)}(\beta^{2}-\|v_0\|^{2})=\|v_n\|_{N L}^{2\cdot2^*_{\mu,s}}-\|v_0\|_{N L}^{2\cdot2^*_{\mu,s}}+o(1),
		\end{equation*}
		and by \eqref{eq5.1} and \cite[Theorem 2]{7thm2} we have,
		\begin{equation}\label{eq5.7}
			\beta^{2(\theta-1)}\lim\limits_{n\to \infty}\|v_n-v_0\|^{2} = d^{2\cdot2^{*}_{\mu ,s}}.
		\end{equation}
		If $d = 0$, then $v_n\to v_0$ in $X_0$ as $\beta>0$. Let, us suppose $d>0$ , by contradiction.
		Reasoning as in Lemma \ref{Lemma4.2}, by \eqref{eq5.1}, \eqref{eq5.7} we get
		\begin{equation}\label{eq5.8}
			\beta^{2}\geq S_H^\frac{2^{*}_{\mu ,s}}{2^{*}_{\mu ,s}-\theta}.
		\end{equation}
		Now, from Lemma \ref{Lemma4.2}, \eqref{eq4.1}, \eqref{eq4.14}, \eqref{eq5.1}, \eqref{eq5.8} and the Young inequality, we have
		\begin{align*}
		D_1-D_2\lambda^\frac{2\theta}{2\theta-1+\gamma}
		 > &J_{n,\lambda}{(v_n)}-\frac{1}{2\cdot{2^*}_{\mu,s}}\langle J_{n,\lambda }^{'}(v_n),(v_n) \rangle\\
			 \geq &\displaystyle\left(\frac{1}{2\theta}-\frac{1}{2\cdot 2_{\mu,s}^*}\right)(\beta^{2\theta}+\|v_0\|^{2\theta})\\
			 &-\lambda\left( \frac{1}{1-\gamma}+\frac{1}{2\cdot 2_{\mu,s}^*}\right)|\Omega|^\frac{2_s^*-1+\gamma}{2_s^*}S_s^\frac{-(1-\gamma)}{2}\|v_0\|^{1-\gamma}\\
			 \geq &\displaystyle \left(\frac{1}{2\theta}-\frac{1}{2\cdot 2_{\mu,s}^*}\right){(S_s^H)}^\frac{2_{\mu,s}^*\theta}{2_{\mu,s}^*-\theta}\\
			&-{\left(\frac{1}{2\theta}-\frac{1}{2\cdot 2_{\mu,s}^*}\right)}^\frac{-(1-\gamma)}{2\theta-1+\gamma}{{\left(\lambda\left( \frac{1}{1-\gamma}+\frac{1}{2\cdot 2_{\mu,s}^*}\right)|\Omega|^\frac{2_s^*-1+\gamma}{2_s^*}S_s^\frac{-(1-\gamma)}{2}\right)}^\frac{2\theta}{2\theta-1+\gamma}}\\
		    = & D_1- \lambda^\frac{2\theta}{2\theta-1+\gamma} D_2,
		\end{align*}
which is a contradiction. Hence, $v_n \to v_0$ in $X_0$ as $n \to \infty$, thus passing limits in \eqref{eq2.5} for $v_n$ we get that $v_0$ satisfies \eqref{eq2.3}. Finally, taking $\phi = v_0^-$ in \eqref{eq2.3} and using \eqref{eq3.14}, we infer that $v_0$ is a positive solution of problem\eqref{eq2.2}. Moreover, $v_0$ is different from $u_0$ we got in Theorem \ref{thm3.2}, since $J_\lambda(v_0) = \lim\limits_{n\to \infty}J_\lambda(v_n) \geq \alpha >0 > J_\lambda(u_0)$. This concludes the proof of Theorem \ref{thm1.1}.\QED

\end{document}